\documentclass[11pt,twoside]{amsart}
\usepackage{amsmath,amssymb,amsthm,mathrsfs,enumitem}
\usepackage[hidelinks,hypertexnames=false,bookmarksnumbered]{hyperref}
\allowdisplaybreaks[3]
\numberwithin{equation}{section}
\newtheorem{theorem}{Theorem}[section]
\newtheorem{proposition}[theorem]{Proposition}
\newtheorem{lemma}[theorem]{Lemma}
\theoremstyle{definition}
\newtheorem{definition}[theorem]{Definition}
\theoremstyle{remark}
\newtheorem{remark}[theorem]{Remark}
\newcommand{\R}{\mathbb R}
\newcommand{\N}{\mathbb N}
\newcommand{\RN}{\mathbb R^N}
\newcommand{\D}{\mathcal D}
\newcommand{\cC}{{\mathcal C}}
\newcommand{\cM}{\mathcal M}
\newcommand{\cS}{{\mathcal S}}
\newcommand{\cE}{\mathcal E}
\newcommand{\dd}{\,dx}
\newcommand{\norm}[2][]{\lVert#2\rVert_{#1}}
\newcommand{\Fone}{\textnormal{(F1)}}
\newcommand{\Ftwo}{\textnormal{(F2)}}
\newcommand{\Fthree}{\textnormal{(F3)}}
\newcommand{\Ffour}{\textnormal{(F4)}}
\newcommand{\cz}{c_{\circ}}
\newcommand{\eps}{\varepsilon}

\title[Sobolev-critical and supercritical Born--Infeld ground states]
{Normalized ground states for the nonlinear Born--Infeld problem with Sobolev-critical
	and supercritical growth}
\author[J. Mederski]{Jaros\l aw Mederski}
\address{J. Mederski, Institute of Mathematics, Polish Academy of Sciences,
	ul. \'Sniadeckich 8, 00-656 Warsaw, Poland;
	and Faculty of Mathematics and Computer Science, Nicolaus Copernicus University,
	ul. Chopina 12/18, 87-100 Toru\'n, Poland}
\email{jmederski@impan.pl}
\author[X. Zeng]{Xiangjian Zeng}
\address{X. Zeng, Institute of Mathematics, Polish Academy of Sciences,
	ul. \'Sniadeckich 8, 00-656 Warsaw, Poland;
	and Faculty of Mathematics, Informatics and Mechanics, University of Warsaw,
	ul. Banacha 2, 02-097 Warsaw, Poland}
\email{xzeng@impan.pl}
\date{}

\begin{document}
	\begin{abstract}
		We study normalized solutions to the nonlinear Born--Infeld problem
		\begin{equation*}
			-{\rm div}\left(\frac{\nabla u}{\sqrt{1-|\nabla u|^2}}\right)
			+\mu u=f(u)
			\quad\hbox{in }\RN,\qquad \int_{\RN}|u|^2\dd=2m^2,
		\end{equation*}
		where $N\ge3$, the mass is prescribed, the frequency $\mu>0$ is
		unknown, and $f$ has Sobolev-critical or supercritical growth at the
		origin, with no upper growth bound at infinity. We find three
		thresholds $c_*\le\rho_0<\widetilde r_e$ for $m^2$. Normalized action
		ground states exist for $m^2>c_*$ and do not exist for $m^2<c_*$. Least
		energy normalized solutions exist for $m^2>\rho_0$, and the energy
		attains its minimum on the mass sphere if and only if
		$m^2\ge\widetilde r_e$. In dimensions $3\le N\le6$, we prove that
		$c_*<\rho_0$ and that the first two existence results hold also for
		$m^2=c_*$ and $m^2=\rho_0$, respectively. The Lagrange multiplier of the
		energy minimizers is obtained from an abstract multiplier rule for
		constrained local minima of functionals with a convex nonsmooth part.
\end{abstract}
	\maketitle
	
	\section{Introduction and main results}\label{sec:intro}
	
	We study the prescribed mass problem
	\begin{equation}\label{eq:normalized}
		\begin{cases}
			-{\rm div}\left(\dfrac{\nabla u}
			{\sqrt{1-|\nabla u|^2}}\right)+\mu u=f(u)&\hbox{in }\RN,\\[5mm]
			\displaystyle\int_{\RN}|u|^2\dd=2m^2,
		\end{cases}
	\end{equation}
	where $N\ge3$, the number $m>0$ is given, and the \emph{frequency} $\mu>0$ is
	unknown and appears as a Lagrange multiplier.
	Let $2^*=2N/(N-2)$. Throughout the paper, unless otherwise stated, we assume that
	\begin{enumerate}[label=\textnormal{(F\arabic*)},leftmargin=2.8em]
		\item\label{ass:F1}
		$f\in \cC^1(\R)$ and $f(0)=f'(0)=0$;
		\item\label{ass:F2}
		the function $f(s)/|s|$ is strictly increasing on $\R\setminus\{0\}$;
		\item\label{ass:F3}
		$\displaystyle\lim_{|s|\to\infty}f(s)/s=+\infty$;
		\item\label{ass:F4}
		for some finite $p\ge2^*$,
		\begin{equation*}
			0<\liminf_{s\to0}\frac{f(s)}{|s|^{p-2}s}
			\le\limsup_{s\to0}\frac{f(s)}{|s|^{p-2}s}<\infty.
		\end{equation*}
	\end{enumerate}
	Assumptions \Fone--\Fthree\ are those of \cite{MZ12}, and they
	correspond to the assumptions of \cite{BIM} in the positive mass case.
	The exponent $p$ in \Ffour\ describes the behavior of $f$ at zero.
	Examples are $f(s)=|s|^{p-2}s$,
	$f(s)=|s|^{p-2}s+|s|^{q-2}s$ with $q>p\ge2^*$, and
	$f(s)=|s|^{p-2}s\,e^{s^2}$.
	We emphasize that $f$ is not assumed to be odd. In particular,
	different positive coefficients on the positive and on the negative
	half-line are allowed.
	
	The differential operator in \eqref{eq:normalized} is the Born--Infeld
	operator. It appears in the electrostatic case of the Born--Infeld
	theory \cite{BI,BdAP}, and it is the mean curvature operator of
	spacelike graphs in the Lorentz--Minkowski space \cite{BaS}.
	The associated functional is
	\begin{equation*}
		\Psi(u)=\int_{\RN}\left(1-\sqrt{1-|\nabla u|^2}\right)\dd,
	\end{equation*}
	which is finite only if $\norm[\infty]{\nabla u}\le1$. For odd $f$,
	solutions of \eqref{eq:normalized} arise formally as profiles of
	standing waves of prescribed $L^2$-norm for the Schr\"odinger and
	Klein--Gordon type evolutions associated with this functional.
	We refer to \cite{MZ12} for these and other motivations.
	
	If the frequency $\mu>0$ is prescribed and the mass is free, the
	equation in \eqref{eq:normalized} was studied by Bieganowski, Ikoma and
	the first author \cite{BIM}. For a large class of nonlinearities, with
	no upper growth restriction at infinity, they obtained \emph{ground states},
	that is, nontrivial solutions with the least action, by minimization on
	a combined Nehari--Poho\v{z}aev set. The regularity theory needed for
	the singular operator is due to Byeon, Ikoma, Malchiodi and Mari
	\cite{BIMM}. The paper \cite{BIM} also treats the zero-frequency
	equation, $\mu=0$, for nonlinearities which are Sobolev-supercritical
	at the origin. For the pure power $|s|^{p-2}s$ with $p>2^*$,
	Luo \cite{Luo} recently proved that the positive ground state of the
	zero-frequency equation is unique up to translations and nondegenerate.
	
	Concerning the prescribed mass problem \eqref{eq:normalized}, Baldelli,
	the first author and Pomponio \cite[Theorem~1.1]{BMP} obtained solutions
	for the pure power $f(s)=|s|^{p-2}s$ in the $L^2$-subcritical range
	$2<p<2+4/N$ by truncating the operator as in \cite{BCD}, and the
	$L^2$-supercritical case was left open there.
	In \cite{MZ12} we considered nonlinearities which behave like
	$|s|^{p-2}s$ at the origin with $2<p<2^*$. The solutions are found among
	the ground states of the prescribed frequency problem: if $a(\mu)$
	denotes the least action at the frequency $\mu$, then the frequency is
	chosen at a local minimum point of
	\begin{equation*}
		b_m(\mu)=a(\mu)-\mu m^2.
	\end{equation*}
	This idea goes back, for semilinear equations, to Dovetta, Serra and
	Tilli \cite{DST}.
	For $2+4/N<p<2^*$ it gives a threshold $c_*>0$ such that ground states
	of mass $2m^2$ exist if and only if $m^2\ge c_*$, and a threshold
	$r_e=\inf_{\mu>0}a(\mu)/\mu>c_*$ such that the energy attains its
	minimum on the mass sphere if and only if $m^2\ge r_e$.
	
	In this paper we consider nonlinearities which are Sobolev-critical or
	Sobolev-supercritical at the origin, that is, $p\ge2^*$, again with no
	upper growth restriction at infinity. For the Laplacian, the corresponding pure
	power problem has no solution: by the Poho\v{z}aev identity, the
	equation $-\Delta u+\mu u=|u|^{p-2}u$ has no nontrivial solution in
	$H^1(\RN)\cap L^p(\RN)$ if $\mu>0$ and $p\ge2^*$, and normalized
	solutions of Sobolev-critical Schr\"odinger equations require a
	subcritical perturbation, see e.g.\ \cite{Soave20crit}. The Born--Infeld
	problem is different because of the gradient constraint. Indeed,
	a function with $\norm[\infty]{\nabla u}\le1$ is bounded in terms of
	its mass, see Lemma \ref{lem:amplitude}. Hence the growth of $f$ at
	infinity plays no role, the energy is bounded from below on every mass
	sphere, and ground states exist at every positive frequency for every
	$p>2$. The difference between $p<2^*$ and $p\ge2^*$ appears at small
	frequencies. For $p<2^*$ the action of the ground states tends to zero
	as $\mu\to0^+$, whereas for $p\ge2^*$ we prove that
	\begin{equation*}
		a_0:=\lim_{\mu\to0^+}a(\mu)>0.
	\end{equation*}
	The positive limit $a_0$ changes the relation between the least action
	at a prescribed frequency and the least energy at a prescribed mass.
	With $F(s)=\int_0^s f(t)\,dt$, the \emph{energy} and the \emph{action}
	at the frequency $\mu>0$ are
	\begin{equation}\label{eq:functionals}
		\cE(u)=\Psi(u)-\int_{\RN}F(u)\dd,\qquad
		I_\mu(u)=\cE(u)+\frac{\mu}{2}\norm[2]{u}^2.
	\end{equation}
	Hence an \emph{action ground state}
	of mass $2m^2$ at the frequency $\mu$ has energy $b_m(\mu)$. Let $e_m$
	be the infimum of $\cE$ on the mass sphere. A function attaining it is
	an \emph{energy ground state}. For $2<p<2^*$ we have $b_m(\mu)\to0$ as
	$\mu\to0^+$, and we proved in \cite{MZ12} that
	$e_m=\inf_{\mu>0}b_m(\mu)$ and that the energy ground states are exactly
	the action ground states whose frequency is a global minimum point of
	$b_m$.
	
	For $p\ge2^*$ we have $b_m(\mu)\to a_0>0$ as $\mu\to0^+$, and the
	picture changes in two respects. First, we show that
	\begin{equation*}
		e_m=\min\Big\{0,\inf_{\mu>0}b_m(\mu)\Big\},
	\end{equation*}
	where the zero is the limit of the energy of $R^{-N/2}\phi(x/R)$,
	$\phi\in\cC_c^\infty(\RN)$, as $R\to\infty$. It cannot be omitted, since
	$\inf_{\mu>0}b_m(\mu)>0=e_m$ for $m^2<\widetilde r_e$, see Theorem
	\ref{thm:least}\ref{it:least-lower}. The threshold $\widetilde r_e$ will
	be defined in Theorem \ref{thm:energy} and characterized by
	\eqref{eq:threshold-minimum-mass}. It coincides with $r_e$ by Remark
	\ref{rem:direct-scope}, and energy ground states exist if and only if
	$m^2\ge\widetilde r_e$, as for $2+4/N<p<2^*$. We obtain them by a direct
	minimization which does not require \Fthree, and their Lagrange
	multiplier by a multiplier rule proved in Appendix
	\ref{app:multiplier}.
	Second, the local minimum points of $b_m$ still give normalized action
	ground states for $m^2>c_*$, where $c_*<\widetilde r_e$. The new feature
	is that $b_m$ may attain a positive global minimum. The corresponding
	ground states do not minimize the energy on the mass sphere, but they
	have the least energy among all solutions of \eqref{eq:normalized}.
	Such \emph{least energy normalized solutions}, in the terminology of
	De Coster, Dovetta, Galant and Serra \cite{DDGS}, exist above a third
	threshold $\rho_0\in[c_*,\widetilde r_e)$, which does not appear for
	$p<2^*$.
	
	\subsection{Variational setting}
	We use the space and the notation of \cite{BIM,MZ12}. Fix $j_0\in\N$ with
	$2j_0>N$, and set
	\begin{equation*}
		X:=\{u\in H^1(\RN):\nabla u\in L^{2j_0}(\RN)\},\qquad
		\norm[X]{u}:=\left(\norm[H^1]{u}^2+
		\norm[2j_0]{\nabla u}^2\right)^{1/2}.
	\end{equation*}
	Then $X\subset H^1(\RN)\cap\cC_0(\RN)$ with continuous inclusion, where
	$\cC_0(\RN)$ is the space of continuous functions vanishing at infinity.
	We regard $\Psi$ as a convex functional on $X$ with values in
	$[0,+\infty]$, setting $\Psi(u)=+\infty$ if $\norm[\infty]{\nabla u}>1$,
	and write
	\begin{equation*}
		\D=\{u\in X:\norm[\infty]{\nabla u}\le1\}.
	\end{equation*}
	The potential term in \eqref{eq:functionals} is finite and of class
	$\cC^1$ on $X$, since on each bounded interval $|F(s)|+|f(s)s|\le C s^2$
	and $|f(s)|\le C|s|$.
	A \emph{critical point} of $I_\mu$ is an element $u\in\D$ satisfying
	the variational inequality of Szulkin \cite{S}:
	\begin{equation}\label{eq:critical}
		\Psi(v)-\Psi(u)+\mu\int_{\RN}u(v-u)\dd
		-\int_{\RN}f(u)(v-u)\dd\ge0
		\qquad\text{for every }v\in X.
	\end{equation}
	By \cite[Lemma~2.3]{BIM}, see also \cite{BIMM}, every critical point $u$
	is a classical solution of the equation in \eqref{eq:normalized}, and it
	satisfies $\norm[\infty]{\nabla u}<1$ and the Poho\v{z}aev identity
	$M_\mu(u)=0$, where
	\begin{equation}\label{eq:M}
		M_\mu(u)=\Psi(u)-\int_{\RN}G(u)\dd
		+\frac{N+2}{2N}\mu\norm[2]{u}^2,
		\qquad G(s)=F(s)+\frac1Nf(s)s.
	\end{equation}
	The least action level and the set of \emph{action ground states} at the
	frequency $\mu$ are
	\begin{equation}\label{eq:a}
		\cM_\mu=\{u\in\D\setminus\{0\}:M_\mu(u)=0\},\qquad
		a(\mu)=\inf_{\cM_\mu}I_\mu,\qquad
		\cC(\mu)=\{u\in\cM_\mu:I_\mu(u)=a(\mu)\}.
	\end{equation}
	In particular, every nontrivial critical point of $I_\mu$ belongs to
	$\cM_\mu$. On the other hand, by Proposition \ref{prop:fixed}, the set
	$\cC(\mu)$ is nonempty for every $\mu>0$ and its elements are critical
	points of $I_\mu$. Hence $\cC(\mu)$ is the set of the nontrivial
	critical points of $I_\mu$ with the least action. We emphasize that no
	symmetry is imposed in its definition. Since neither the uniqueness of
	the action ground states nor the differentiability of $a$ is known,
	we consider, as in \cite{CGT,MZ12}, the extreme masses in $\cC(\mu)$.
	For $m>0$ and $\mu>0$ we set
	\begin{align}
		\cS_m&=\{u\in\D:\tfrac12\norm[2]{u}^2=m^2\},&
		e_m&=\inf_{u\in\cS_m}\cE(u),\label{eq:mass}\\
		c_-(\mu)&=\inf_{u\in\cC(\mu)}\tfrac12\norm[2]{u}^2,&
		c_+(\mu)&=\sup_{u\in\cC(\mu)}\tfrac12\norm[2]{u}^2,\label{eq:cp}\\
		b_m(\mu)&=a(\mu)-\mu m^2,&
		B_m&=\inf_{\mu>0}b_m(\mu).\label{eq:b}
	\end{align}
	Thus all the thresholds below refer to half of the squared $L^2$-norm.
	
	\begin{definition}\label{def:solutions}
		Let $m>0$.
		\begin{enumerate}[label=\textnormal{(\alph*)}]
			\item\label{defitem:sol} A \emph{solution} of \eqref{eq:normalized} is a pair
			$(\mu,u)\in(0,\infty)\times\cS_m$ such that $u$ is a critical point of
			$I_\mu$.
			\item\label{defitem:action} A \emph{normalized action ground state} is an element of
			$\cC(\mu)\cap\cS_m$ for some $\mu>0$.
			\item\label{defitem:energy} An \emph{energy ground state} is a function $u\in\cS_m$ with
			$\cE(u)=e_m$.
			\item\label{defitem:least} A \emph{least energy normalized solution} is a solution
			$(\mu,u)$ of \eqref{eq:normalized} such that $\cE(u)\le\cE(v)$ for
			every solution $(\nu,v)$ of \eqref{eq:normalized}.
		\end{enumerate}
	\end{definition}
	If $u\in\cC(\mu)\cap\cS_m$ is a normalized action ground state, then
	$(\mu,u)$ is a solution, because every element of $\cC(\mu)$ is a
	critical point of $I_\mu$, and $\cE(u)=b_m(\mu)$.

	The zero-frequency problem enters through
	\begin{equation}\label{eq:M0}
		M_0(u)=\Psi(u)-\int_{\RN}G(u)\dd,\qquad
		\cM_0=\{u\in\D\setminus\{0\}:M_0(u)=0\},
	\end{equation}
	which are the counterparts of \eqref{eq:M} and \eqref{eq:a} for $\mu=0$.
	
	For some results we use the additional condition
	\begin{equation}\label{eq:critical-comparison}
		sf(s)\le 2^*F(s)\qquad\text{for every }s\in\R.
	\end{equation}
	The critical pure power satisfies it with equality, and another
	example is given in Remark \ref{rem:critical-example}. Note that
	\eqref{eq:critical-comparison} forces $p=2^*$ in \Ffour. Indeed,
	it implies that $F(s)/|s|^{2^*}$ is nonincreasing in $|s|$ on each of
	the half-lines $(-\infty,0)$ and $(0,\infty)$, so that
	$F(s)\ge\min\{F(1),F(-1)\}|s|^{2^*}$ for $|s|\le1$, whereas
	$F(s)=O(|s|^{p})$ as $s\to0$ by \Ffour.
	
	\subsection{Main results}
	Our first result concerns the nonexistence of solutions for small masses
	and the existence of normalized action ground states.
	
	\begin{theorem}\label{thm:action}
		Assume \Fone--\Ffour, and set
		\begin{equation}\label{eq:thresholds}
			\cz:=\inf_{\mu>0}\,\inf_{u\in\cM_\mu}\tfrac12\norm[2]{u}^2,\qquad
			c_*:=\inf_{\mu>0}c_-(\mu).
		\end{equation}
		Then the limit $a_0:=\lim_{\mu\to0^+}a(\mu)$ exists, and
		\begin{equation}\label{eq:mainthresholds}
			0<a_0<\infty,\qquad 0<\cz<c_*<\infty.
		\end{equation}
		Moreover, for $m>0$ the following assertions hold.
		\begin{enumerate}[label=\textnormal{(\roman*)}]
			\item\label{it:none}
			If $m^2\le\cz$, then \eqref{eq:normalized} has no solution.
			\item\label{it:below}
			If $m^2<c_*$, then there is no normalized action ground state, that
			is, $\cC(\mu)\cap\cS_m=\varnothing$ for every $\mu>0$.
			\item\label{it:above}
			If $m^2>c_*$, then there is $\mu>0$ such that
			$c_-(\mu)=c_+(\mu)=m^2$, that is, $\cC(\mu)\subset\cS_m$.
			\item\label{it:at}
			If $m^2=c_*$, then a normalized action ground state exists if and only
			if the infimum defining $c_*$ is attained. Moreover, if
			$c_-(\mu_0)=c_*$ for some $\mu_0>0$, then
			$\cC(\mu_0)\cap\cS_m\ne\varnothing$, and every $u$ in this set satisfies
			$\cE(u)=a(\mu_0)-\mu_0c_*\ge a_0$.
		\end{enumerate}
		Furthermore, for every $\mu>0$, every $u\in\cC(\mu)$ has the form
		$u(x)=\pm W(|x-x_0|)$ with some $x_0\in\RN$, where $W$ is positive and strictly
		decreasing on $[0,\infty)$. If $f$ is odd, then the normalized action ground states in
		\ref{it:above} and \ref{it:at} can be chosen positive.
	\end{theorem}
	
	The nonexistence in
	\ref{it:none} also holds for $2+4/N\le p<2^*$, as shown in Remark
	\ref{rem:smallmass}, and the attainment of $c_*$ is discussed in
	Theorem \ref{thm:attainment} below.
	
	The second result concerns the minimization of the energy on the mass
	sphere.
	
	\begin{theorem}[Energy ground states]\label{thm:energy}
		Assume \Fone, \Ftwo\ and \Ffour, and set
		\begin{equation*}
			\widetilde r_e:=\inf\{m^2>0:e_m<0\}.
		\end{equation*}
		Then $0<\widetilde r_e<\infty$, and
		\begin{equation}\label{eq:threshold-minimum-mass}
			\widetilde r_e=\min\left\{\tfrac12\norm[2]{v}^2:
			v\in\D\setminus\{0\},\ \cE(v)\le0\right\}.
		\end{equation}
		Moreover, for $m>0$ the following assertions hold.
		\begin{enumerate}[label=\textnormal{(\roman*)}]
			\item\label{it:en-below}
			If $m^2<\widetilde r_e$, then $e_m=0$ is not attained.
			\item\label{it:en-at}
			If $m^2=\widetilde r_e$, then $e_m=0$ is attained.
			\item\label{it:en-above}
			If $m^2>\widetilde r_e$, then $e_m<0$ is attained.
		\end{enumerate}
		Furthermore, every energy ground state $u\in\cS_m$ has constant sign,
		and $(\mu,u)$ is a solution of \eqref{eq:normalized} for some $\mu>0$.
		Whenever $e_m$ is attained, there is an energy ground state of the
		form $u(x)=\pm W(|x|)$, where $W$ is positive and nonincreasing on
		$[0,\infty)$.
	\end{theorem}
	
	If \Fthree\ holds as well, then for $m^2\ge\widetilde r_e$ the energy
	ground states are exactly the action ground states whose frequency is a
	global minimum point of $b_m$, see Proposition \ref{prop:duality}. In
	particular, they have the form described in Theorem \ref{thm:action}.
	The next theorem shows that the global minimum points of $b_m$ keep a
	variational meaning also for $m^2<\widetilde r_e$.
	
	\begin{theorem}[Least energy normalized solutions]\label{thm:least}
		Assume \Fone--\Ffour, let $a_0$, $c_*$ and $\widetilde r_e$ be given by
		Theorems \ref{thm:action} and \ref{thm:energy}, and set
		\begin{equation*}
			\rho_0:=\inf_{\mu>0}\frac{a(\mu)-a_0}{\mu}.
		\end{equation*}
		Then $c_*\le\rho_0<\widetilde r_e$. Moreover, for $m>0$ the following
		assertions hold.
		\begin{enumerate}[label=\textnormal{(\roman*)}]
			\item\label{it:least-lower}
			Every solution $(\mu,u)$ of \eqref{eq:normalized} satisfies
			$\cE(u)\ge b_m(\mu)\ge B_m$. Moreover, $B_m>0$ if
			$m^2<\widetilde r_e$, and $\cE(u)>a_0$ if $m^2<\rho_0$.
			\item\label{it:least-attain}
			The function $b_m$ attains its infimum on $(0,\infty)$ if and only if
			$m^2>\rho_0$, or $m^2=\rho_0$ and the infimum defining $\rho_0$ is
			attained. In this case $B_m$ is the least energy of the solutions of
			\eqref{eq:normalized}, every global minimum point $\mu$ of $b_m$
			satisfies $\cC(\mu)\subset\cS_m$, and a solution $(\mu,u)$
			of \eqref{eq:normalized} is a least energy normalized solution if and
			only if $b_m(\mu)=B_m$ and $u\in\cC(\mu)$.
		\end{enumerate}
	\end{theorem}
	
	Thus, for $\rho_0<m^2<\widetilde r_e$ the least energy normalized
	solutions are action ground states of positive energy, whereas $e_m=0$
	is not attained, and for $m^2\ge\widetilde r_e$ they coincide with the
	energy ground states.
	
	It remains to determine whether $c_*$ and $\rho_0$ are attained at
	positive frequencies and whether $c_*<\rho_0$. A minimizing sequence
	of frequencies for either threshold is bounded above. To rule out
	convergence to zero, we study the limiting equation
	\begin{equation}\label{eq:zerofreq}
		-{\rm div}\left(\frac{\nabla w}{\sqrt{1-|\nabla w|^2}}\right)=f(w)
		\qquad\hbox{in }\RN.
	\end{equation}
	\begin{theorem}[Attainment of the mass thresholds]\label{thm:attainment}
		Assume \Fone--\Ffour.
		\begin{enumerate}[label=\textnormal{(\roman*)}]
			\item\label{it:attain}
			If $3\le N\le6$, or if \eqref{eq:critical-comparison} holds, then the
			infima defining $c_*$ in \eqref{eq:thresholds} and $\rho_0$ in Theorem
			\ref{thm:least} are attained, and $c_*<\rho_0$.
			\item\label{it:notattain}
			If the infimum defining $c_*$ is not attained, then $N\ge7$ and the
			zero-frequency equation
			\eqref{eq:zerofreq} has a classical solution $w\in\cM_0$, radial and of
			constant sign, such that $\tfrac12\norm[2]{w}^2=c_*$ and $\cE(w)=a_0$.
		\end{enumerate}
	\end{theorem}
	
	The function $w$ in \ref{it:notattain} minimizes $\cE$ on $\cM_0$ by
	\eqref{eq:a0}. Thus, if no action ground state at a positive frequency
	has mass $2c_*$, then this mass is realized by a ground state at zero
	frequency. In dimensions $3\le N\le6$, Theorems
	\ref{thm:action}--\ref{thm:attainment} give the strict ordering
	\begin{equation*}
		0<\cz<c_*<\rho_0<\widetilde r_e.
	\end{equation*}
	Moreover, in these dimensions normalized action ground states exist if
	and only if $m^2\ge c_*$, and the function $b_m$ attains its infimum if
	and only if $m^2\ge\rho_0$.
	In general we only know that $c_*\le\rho_0$. 
	
	\begin{remark}\label{rem:critical-example}
		Condition \eqref{eq:critical-comparison} is satisfied by nonlinearities
		other than the critical pure power. For example, let
		\begin{equation*}
			f(s)=\frac{|s|^{2^*-2}s}{1+|s|^\alpha},
			\qquad 0<\alpha<2^*-2.
		\end{equation*}
		Then \Fone--\Ffour\ hold with $p=2^*$. Indeed, $f$ is odd and of class
		$\cC^1$, with $f(0)=f'(0)=0$,
		$f(s)/(|s|^{2^*-2}s)\to1$ as $s\to0$, and
		$f(s)/s\sim|s|^{2^*-2-\alpha}\to+\infty$ as $|s|\to\infty$.
		Moreover, for $s>0$ the derivative of $f(s)/s$ is
		\begin{equation*}
			\frac{(2^*-2)s^{2^*-3}+(2^*-2-\alpha)s^{2^*+\alpha-3}}
			{(1+s^\alpha)^2}>0,
		\end{equation*}
		which gives \Ftwo. Finally, if $0<t<s$, then
		\begin{equation*}
			f(t)\ge f(s)\left(\frac ts\right)^{2^*-1},
			\qquad
			F(s)=\int_0^s f(t)\,dt\ge\frac{sf(s)}{2^*}.
		\end{equation*}
		Since $F$ and $sf(s)$ are even, this proves
		\eqref{eq:critical-comparison} on $\R$.
		By Proposition \ref{prop:divergence} and Theorem
		\ref{thm:attainment}, in every dimension $N\ge3$ we have
		$c_-(\mu)\to+\infty$ as $\mu\to0^+$, the infima defining $c_*$ and
		$\rho_0$ are attained, and $c_*<\rho_0$. Thus the masses of the ground
		states are unbounded as $\mu\to0^+$ not only for the critical pure
		power.
	\end{remark}
	
	\begin{remark}\label{rem:direct-scope}
		In \cite[Theorem~1.3(iii)]{MZ12}, where $2+4/N\le p<2^*$, the threshold
		for the energy ground states is $r_e=\inf_{\mu>0}a(\mu)/\mu$, which is
		defined by means of the action ground state level, whereas the
		definition of $\widetilde r_e$ in Theorem \ref{thm:energy} involves only
		the levels $e_m$. If \Fone--\Fthree\ hold, then $\widetilde r_e=r_e$.
		Indeed, by Proposition \ref{prop:duality}, we have $e_m<0$ if and only
		if $a(\mu)-\mu m^2<0$ for some $\mu>0$, that is, if and only if
		$m^2>r_e$. In particular, under \Fone--\Ffour\ the number $r_e$ is
		characterized by \eqref{eq:threshold-minimum-mass}, and the threshold
		$\rho_0$ of Theorem \ref{thm:least} would coincide with $\widetilde r_e$
		if $a_0$ were zero.
	\end{remark}
	
	\subsection{Comments and related results}
	\emph{Comparison with the Sobolev-subcritical case.}
	For $2+4/N<p<2^*$, by \cite[Lemma~3.2]{MZ12}, there is $c>0$ such that
	$c_-(\mu)\ge c\mu^{\beta}$ for small $\mu>0$, where
	$\beta=\frac{2}{p-2}-\frac N2<0$. Hence $c_-(\mu)\to\infty$ as
	$\mu\to0^+$, and the infimum $c_*$ is attained. Such a bound cannot hold
	for $p\ge2^*$. Indeed, in this case $\beta\le-1$, whereas $c_-$ is
	integrable near zero by \eqref{eq:integral0}. For $p\ge2^*$ the
	divergence $c_-(\mu)\to\infty$ as $\mu\to0^+$ still holds if
	$N\in\{3,4\}$ or if \eqref{eq:critical-comparison} holds, see
	Proposition \ref{prop:divergence}. For $N\ge5$, Proposition
	\ref{prop:bounded-mass-limit} describes the limits of ground states
	whose masses stay bounded as $\mu\to0^+$.
	
	\emph{Semilinear analogues.}
	A positive limit of the action level at zero frequency also occurs
	for semilinear equations: \cite[Proposition~3.3(ii)]{CGT} establishes
	it for nonlinearities which are Sobolev-critical or supercritical at
	zero and subcritical at infinity. The coexistence of energy
	minimizers above a mass threshold with normalized solutions of
	positive energy below it is known for the Schr\"odinger equation with
	a focusing cubic and a defocusing quintic nonlinearity
	\cite{KOPV,JL22,CGT}, and for the $(2,q)$-Laplacian between the two
	mass-critical exponents \cite{BaIk}, where the dimension also decides
	whether zero-mass solutions have finite $L^2$-norm. In the
	Born--Infeld problem these phenomena are caused by the gradient
	constraint alone, for a single focusing nonlinearity. 
	
	\subsection{Open problems}
	We list some questions which are left open.
	\begin{enumerate}[label=\textnormal{(\arabic*)},leftmargin=2em]
		\item Does \eqref{eq:normalized} have solutions, necessarily
		different from action ground states, for $\cz<m^2<c_*$?
		\item Does \eqref{eq:normalized} have a least energy normalized
		solution for $c_*<m^2<\rho_0$, where the function $b_m$ does not
		attain its infimum?
		\item Suppose that $c_-(\mu_1)<m^2<\liminf_{\mu\to0^+}c_-(\mu)$ for some
		$\mu_1>0$. Then $b_m$ has a local maximum at some
		$\mu_0\in(0,\mu_1]$. Does $\cC(\mu_0)$ contain a state of mass
		$2m^2$, giving a second solution at a different frequency?
	\end{enumerate}
	
	\subsection{Organization of the paper}
	The paper is organized as follows. In Section \ref{sec:tools}, we
	recall the results of \cite{BIM,MZ12} on the prescribed frequency
	problem. Section
	\ref{sec:small} is devoted to the zero-frequency level and to the
	thresholds, and proves Theorem \ref{thm:action} except for the
	inequality $\cz<c_*$. In Section \ref{sec:energy}, we
	prove Theorem \ref{thm:energy} by a direct minimization on the mass
	sphere, relate the energy ground states to the action ground states,
	and prove Theorem \ref{thm:least}. In Section \ref{sec:endpoints}, we
	study limits of ground states as $\mu\to0^+$ and prove both
	$\cz<c_*$ and Theorem \ref{thm:attainment}. Finally, Appendix \ref{app:multiplier}
	contains a multiplier rule for constrained local minima of functionals
	with a convex part, which is used in Section \ref{sec:energy}.
	
	\section{Preliminaries}\label{sec:tools}
	
	Throughout the paper, $C$ and $c$ denote positive constants which may
	change from line to line, and $\omega_N$ is the measure of the unit
	sphere in $\RN$. We set
	\begin{equation*}
		K(s)=\frac12 f(s)s-F(s),\qquad s\in\R,
	\end{equation*}
	and for $u\in\D$ and $\theta>0$ we write
	$u_\theta(x)=\theta u(x/\theta)$. Then $u_\theta\in\D$, and
	\begin{equation}\label{eq:scaling}
		\Psi(u_\theta)=\theta^N\Psi(u),\qquad
		\norm[2]{u_\theta}^2=\theta^{N+2}\norm[2]{u}^2,\qquad
		\int_{\RN}F(u_\theta)\dd=\theta^N\int_{\RN}F(\theta u)\dd,
	\end{equation}
	and the same holds for $G$ and $K$ in place of $F$. We first collect
	some elementary consequences of the assumptions.
	
	\begin{lemma}\label{lem:f}
		Assume \Fone--\Ftwo.
		\begin{enumerate}[label=\textnormal{(\alph*)}]
			\item\label{it:f-sign} For $s\ne0$ we have $f(s)s>0$,
			$0<F(s)<\frac12f(s)s$ and $K(s)>0$.
			\item\label{it:f-convex} $f$ is increasing on $\R$, hence $F$ is convex.
			\item\label{it:f-mono} For every $s\ne0$ the functions
			$\theta\mapsto\theta^{-2}F(\theta s)$,
			$\theta\mapsto\theta^{-2}G(\theta s)$ and $\theta\mapsto K(\theta s)$ are
			nondecreasing on $(0,\infty)$.
			\item\label{it:f-local} For every $R>0$ there is $C>0$ such that
			$F(s)\le G(s)\le f(s)s\le Cs^2$ for $|s|\le R$. If \Ffour\ holds, then
			moreover $f(s)s\le C|s|^{q}$ for $|s|\le R$ and every $q\in[2,p]$.
		\end{enumerate}
	\end{lemma}
	\begin{proof}
		By \Fone, $f(s)/s\to0$ as $s\to0$. Together with \Ftwo, this shows that
		$f(s)/s$ is positive for $s\ne0$, increasing on $(0,\infty)$ and
		decreasing on $(-\infty,0)$. Since
		\begin{equation*}
			F(s)=\int_0^s\frac{f(t)}{t}\,t\,dt,
		\end{equation*}
		we get \ref{it:f-sign}. Moreover, $f(s)=\frac{f(s)}{s}\,s$ is increasing
		on each of the half-lines $(-\infty,0)$ and $(0,\infty)$, and
		$f(s)<0<f(t)$ for $s<0<t$, which is \ref{it:f-convex}. By a change of
		variables,
		\begin{equation*}
			\theta^{-2}F(\theta s)=\int_0^s\frac{f(\theta t)}{\theta t}\,t\,dt,\qquad
			\theta^{-2}G(\theta s)=\theta^{-2}F(\theta s)
			+\frac1N\,\frac{f(\theta s)}{\theta s}\,s^2,
		\end{equation*}
		and both expressions are nondecreasing in $\theta$. Furthermore,
		\begin{equation*}
			sK'(s)=\frac s2\big(f'(s)s-f(s)\big)
			=\frac{s^3}{2}\,\frac{d}{ds}\left(\frac{f(s)}{s}\right)\ge0
			\qquad\text{for }s\ne0,
		\end{equation*}
		so that $K$ is nondecreasing on $(0,\infty)$ and nonincreasing on
		$(-\infty,0)$, and \ref{it:f-mono} follows. Finally,
		\ref{it:f-sign} gives $F(s)\le G(s)\le f(s)s$, and $f(s)s\le Cs^2$ for
		$|s|\le R$ since $f\in\cC^1(\R)$ and $f(0)=0$. Under \Ffour\ there are
		$C,s_0>0$ such that $f(s)s\le C|s|^{p}$ for $|s|\le s_0$, and the
		estimate extends to $|s|\le R$ and to the exponents $q\in[2,p]$ by
		enlarging the constant.
	\end{proof}
	
	Since
	$I_\mu=\frac{2}{N+2}\big(\Psi+\int_{\RN}K(\cdot)\dd\big)+\frac{N}{N+2}M_\mu$
	on $\D$, by Lemma \ref{lem:f}\ref{it:f-sign} we have
	\begin{equation}\label{eq:energyidentity}
		I_\mu(u)=\frac{2}{N+2}
		\left(\Psi(u)+\int_{\RN}K(u)\dd\right)\ge\frac{2}{N+2}\Psi(u),
		\qquad u\in\cM_\mu,
	\end{equation}
	and the same identity holds for $\cE$ on $\cM_0$.
	
	\begin{lemma}\label{lem:amplitude}
		There is $C>0$, depending only on $N$, such that every $u\in\D$ satisfies
		\begin{equation}\label{eq:amplitude}
			\norm[\infty]{u}^{N+2}\le C\norm[2]{u}^2,
			\qquad
			\norm[\infty]{u}^{N}\le C\norm[2]{\nabla u}^{2},
		\end{equation}
		and
		\begin{equation}\label{eq:kinetic}
			\frac12\norm[2]{\nabla u}^2\le\Psi(u)\le\norm[2]{\nabla u}^2.
		\end{equation}
	\end{lemma}
	\begin{proof}
		The estimates \eqref{eq:amplitude} are proved in
		\cite[Lemma~A.1]{MZ12}, and \eqref{eq:kinetic} follows from
		$1-\sqrt{1-t^2}=t^2/(1+\sqrt{1-t^2})$ for $0\le t\le1$.
	\end{proof}
	
	The following proposition collects the results on the prescribed
	frequency problem which we use, see
	\cite[Theorem~1.6, Lemmas~5.2 and~5.5]{BIM} and
	\cite[Subsection~1.1, Section~2, Lemma~3.1 and Theorem~B.5]{MZ12}, where
	only \Fone--\Fthree\ are used.
	
	\begin{proposition}\label{prop:fixed}
		Assume \Fone--\Fthree.
		\begin{enumerate}[label=\textnormal{(\roman*)}]
			\item\label{it:exist}
			For every $\mu>0$, $a(\mu)>0$ and $\cC(\mu)\ne\varnothing$. Every
			$u\in\cC(\mu)$ is a critical point of $I_\mu$, and it has the form
			$u(x)=\pm W(|x-x_0|)$ with some $x_0\in\RN$, where
			$W\in\cC^2([0,\infty))$, $W>0$ on $[0,\infty)$ and $W'<0$ on
			$(0,\infty)$.
			\item\label{it:fibres}
			Let $\mu>0$ and $u\in\D\setminus\{0\}$. Then there is a unique
			$\theta_\mu(u)>0$ such that $u_{\theta_\mu(u)}\in\cM_\mu$. Moreover,
			$M_\mu(u_\theta)>0$ for $0<\theta<\theta_\mu(u)$,
			$M_\mu(u_\theta)<0$ for $\theta>\theta_\mu(u)$, and
			$I_\mu(u_\theta)<I_\mu(u_{\theta_\mu(u)})$ for
			$\theta\ne\theta_\mu(u)$.
			\item\label{it:compact}
			If $\mu_n\to\mu>0$ and $u_n\in\cC(\mu_n)$, then there are $y_n\in\RN$
			and $u\in\cC(\mu)$ such that, up to a subsequence,
			$u_n(\cdot+y_n)\to u$ in $X$. In particular, the infimum and the
			supremum in \eqref{eq:cp} are attained for every $\mu>0$, and on
			$(0,\infty)$ the function $c_-$ is lower semicontinuous and $c_+$ is
			upper semicontinuous.
			\item\label{it:derivatives}
			The function $a$ is continuous and nondecreasing on $(0,\infty)$, and
			at every $\mu>0$ it has the one-sided derivatives
			\begin{equation}\label{eq:derivatives}
				a'_+(\mu)=c_-(\mu),\qquad a'_-(\mu)=c_+(\mu).
			\end{equation}
			\item\label{it:large}
			$c_-(\mu)\to+\infty$ and $a(\mu)/\mu\to+\infty$ as $\mu\to\infty$.
			\item\label{it:selection}
			Let $m>0$. By \eqref{eq:derivatives}, every local minimum point $\mu_0$
			of $b_m$ satisfies $c_-(\mu_0)=c_+(\mu_0)=m^2$, that is,
			$\cC(\mu_0)\subset\cS_m$.
		\end{enumerate}
	\end{proposition}
	
	\section{Small frequencies and the thresholds}\label{sec:small}
	
	Recall the definition \eqref{eq:M0} of $M_0$ and $\cM_0$.
	
	\begin{lemma}\label{prop:zero-level}
		Assume \Fone--\Fthree.
		\begin{enumerate}[label=\textnormal{(\alph*)}]
			\item\label{it:zl-fibre}
			For every $u\in\D\setminus\{0\}$ there is a unique $\theta_0(u)>0$ such
			that $u_{\theta_0(u)}\in\cM_0$. Moreover, $M_0(u_\theta)>0$ for
			$0<\theta<\theta_0(u)$, $M_0(u_\theta)<0$ for $\theta>\theta_0(u)$, and
			$\cE(u_\theta)<\cE(u_{\theta_0(u)})$ for $\theta\ne\theta_0(u)$.
			\item\label{it:zl-limit}
			$a(\mu)\ge\inf_{\cM_0}\cE$ for every $\mu>0$, the limit
			$a_0:=\lim_{\mu\to0^+}a(\mu)$ exists, and
			\begin{equation}\label{eq:a0}
				a_0=\inf_{\cM_0}\cE.
			\end{equation}
		\end{enumerate}
	\end{lemma}
	\begin{proof}
		\ref{it:zl-fibre} By \eqref{eq:scaling},
		\begin{equation*}
			\theta^{-N}M_0(u_\theta)=\Psi(u)-\int_{\RN}G(\theta u)\dd,
			\qquad
			\frac{d}{d\theta}\cE(u_\theta)
			=N\theta^{N-1}\Big(\Psi(u)-\int_{\RN}G(\theta u)\dd\Big).
		\end{equation*}
		Since $u$ is bounded, Lemma \ref{lem:f}\ref{it:f-local} and the
		dominated convergence theorem show that
		$\theta\mapsto\int_{\RN}G(\theta u)\dd$ is continuous. By Lemma
		\ref{lem:f}\ref{it:f-mono}, the function
		$\theta\mapsto\theta^{-2}\int_{\RN}G(\theta u)\dd$ is positive and
		nondecreasing. Hence $\int_{\RN}G(\theta u)\dd$ is increasing in
		$\theta$, it is not larger than $\theta^2\int_{\RN}G(u)\dd$ for
		$\theta\le1$, and it is not smaller than $\theta^2\int_{\RN}G(u)\dd$
		for $\theta\ge1$. Since $\Psi(u)>0$, the assertion follows.
		
		\ref{it:zl-limit} Let $\mu>0$ and $u\in\cC(\mu)$. Then
		$M_0(u)=-\frac{N+2}{2N}\mu\norm[2]{u}^2<0$, hence
		$\theta:=\theta_0(u)<1$ by \ref{it:zl-fibre}. By
		\eqref{eq:energyidentity}, \eqref{eq:scaling} and Lemma
		\ref{lem:f}\ref{it:f-mono},
		\begin{equation*}
			\inf_{\cM_0}\cE\le\cE(u_{\theta})
			=\frac{2\theta^N}{N+2}\Big(\Psi(u)+\int_{\RN}K(\theta u)\dd\Big)
			\le\frac{2}{N+2}\Big(\Psi(u)+\int_{\RN}K(u)\dd\Big)=a(\mu).
		\end{equation*}
		Conversely, let $u\in\cM_0$ and $0<\mu\le1$. Then
		$M_\mu(u)=\frac{N+2}{2N}\mu\norm[2]{u}^2>0$, and
		$M_\mu(u_\theta)\le M_1(u_\theta)<0$ for $\theta>\theta_1(u)$. By
		Proposition \ref{prop:fixed}\ref{it:fibres},
		$1<\theta_\mu(u)\le\theta_1(u)$. Since $\cE(u_\theta)\le\cE(u)$ for
		every $\theta>0$ by \ref{it:zl-fibre}, we obtain
		\begin{equation*}
			a(\mu)\le I_\mu(u_{\theta_\mu(u)})
			=\cE(u_{\theta_\mu(u)})+\frac\mu2\theta_\mu(u)^{N+2}\norm[2]{u}^2
			\le\cE(u)+\frac\mu2\theta_1(u)^{N+2}\norm[2]{u}^2.
		\end{equation*}
		Letting $\mu\to0^+$ and taking the infimum over $u\in\cM_0$, we get
		$\lim_{\mu\to0^+}a(\mu)\le\inf_{\cM_0}\cE$, where the limit exists
		because $a$ is nondecreasing.
	\end{proof}
	
	\begin{lemma}\label{prop:positive}
		Assume \Fone--\Ffour, and let $a_0$ and $\cz$ be given by \eqref{eq:a0}
		and \eqref{eq:thresholds}. Then $a_0>0$, and
		\begin{equation}\label{eq:smallmass-characterization}
			\cz=\inf_{u\in\cM_0}\frac12\norm[2]{u}^2>0.
		\end{equation}
		Moreover, for every $\mu>0$ and every $u\in\cM_\mu$,
		\begin{equation}\label{eq:smallmass-strict}
			\frac12\norm[2]{u}^2>\cz.
		\end{equation}
	\end{lemma}
	\begin{proof}
		Let $u\in\cM_0$ with $\Psi(u)\le1$. By Lemma \ref{lem:amplitude},
		$\norm[\infty]{u}^N\le C\norm[2]{\nabla u}^2\le C\Psi(u)\le C$. Hence
		Lemma \ref{lem:f}\ref{it:f-local} with $q=2^*$, the Sobolev inequality
		and \eqref{eq:kinetic} give
		\begin{equation*}
			\Psi(u)=\int_{\RN}G(u)\dd\le C\norm[2^*]{u}^{2^*}
			\le C\norm[2]{\nabla u}^{2^*}\le C\,\Psi(u)^{2^*/2}.
		\end{equation*}
		Since $\Psi(u)>0$ and $2^*>2$, it follows that $\Psi$ is bounded away
		from zero on $\cM_0$, and \eqref{eq:energyidentity} and \eqref{eq:a0}
		yield $a_0>0$.
		
		Let now $\mu>0$ and $u\in\cM_\mu$ with $\norm[2]{u}\le1$. Then
		$\norm[\infty]{u}\le C$ by Lemma \ref{lem:amplitude}, and $M_\mu(u)=0$,
		\eqref{eq:kinetic}, Lemma \ref{lem:f}\ref{it:f-local} with $q=2+4/N$
		and the Gagliardo--Nirenberg inequality give
		\begin{equation*}
			\frac12\norm[2]{\nabla u}^2
			\le\Psi(u)
			<\int_{\RN}G(u)\dd
			\le C\norm[2+4/N]{u}^{2+4/N}
			\le C\norm[2]{\nabla u}^2\norm[2]{u}^{4/N}.
		\end{equation*}
		Since $\nabla u\ne0$, we get $\norm[2]{u}^{4/N}\ge c$. Thus $\cz>0$.
		
		It remains to prove the identity in
		\eqref{eq:smallmass-characterization} and the inequality
		\eqref{eq:smallmass-strict}.
		Let $\mu>0$ and $u\in\cM_\mu$. Then
		$M_0(u)=-\frac{N+2}{2N}\mu\norm[2]{u}^2<0$, so that $\theta_0(u)<1$ by
		Lemma \ref{prop:zero-level}\ref{it:zl-fibre}, and
		\begin{equation}\label{eq:smallmass-proof}
			\inf_{v\in\cM_0}\frac12\norm[2]{v}^2
			\le\frac12\norm[2]{u_{\theta_0(u)}}^2
			=\theta_0(u)^{N+2}\,\frac12\norm[2]{u}^2<\frac12\norm[2]{u}^2 .
		\end{equation}
		Taking the infimum over $\mu>0$ and $u\in\cM_\mu$ we get
		$\inf_{v\in\cM_0}\frac12\norm[2]{v}^2\le\cz$. Conversely, let
		$v\in\cM_0$ and $\theta>1$. Then $M_0(v_\theta)<0$ by Lemma
		\ref{prop:zero-level}\ref{it:zl-fibre}, hence $v_\theta\in\cM_\mu$ for
		$\mu=-2NM_0(v_\theta)/\big((N+2)\norm[2]{v_\theta}^2\big)>0$, and
		$\cz\le\frac12\norm[2]{v_\theta}^2=\theta^{N+2}\frac12\norm[2]{v}^2$.
		Letting $\theta\to1^+$ we get $\cz\le\frac12\norm[2]{v}^2$. This proves
		\eqref{eq:smallmass-characterization}, and then
		\eqref{eq:smallmass-proof} gives \eqref{eq:smallmass-strict}.
	\end{proof}
	
	\begin{remark}\label{rem:smallmass}
		The identity in \eqref{eq:smallmass-characterization} and the inequality
		\eqref{eq:smallmass-strict} use only \Fone--\Fthree, and the proof of
		$\cz>0$ uses \Ffour\ only through the bound $G(s)\le C|s|^{2+4/N}$ for
		small $|s|$. Hence $\cz>0$ whenever \Fone--\Fthree\ hold and
		$f(s)=O(|s|^{1+4/N})$ as $s\to0$, in particular in the setting of
		\cite{MZ12} with $2+4/N\le p<2^*$. Since every solution $(\mu,u)$ of
		\eqref{eq:normalized} satisfies $u\in\cM_\mu$, problem
		\eqref{eq:normalized} has no solution with $m^2\le\cz$ in this case
		either. This complements the nonexistence result for normalized action
		ground states in \cite[Theorem~1.1]{MZ12}.
	\end{remark}

	Recall that $\cz$ and $c_*$ are defined in \eqref{eq:thresholds}, and
	$\rho_0$ in Theorem \ref{thm:least}. For $\mu>0$ we write
	\begin{equation}\label{eq:q}
		\rho(\mu):=\frac{a(\mu)-a_0}{\mu},\qquad\text{so that}\qquad
		b_m(\mu)=a_0+\mu\big(\rho(\mu)-m^2\big)
		\quad\text{and}\quad\rho_0=\inf_{\mu>0}\rho(\mu).
	\end{equation}
	
	\begin{lemma}\label{prop:thresholds}
		Assume \Fone--\Ffour. Then the function $c_-$ is integrable on
		$(0,\mu)$ for every $\mu>0$, and
		\begin{equation}\label{eq:integral0}
			a(\mu)=a_0+\int_0^\mu c_-(t)\,dt\qquad\text{for every }\mu>0 .
		\end{equation}
		In particular,
		\begin{equation}\label{eq:gap}
			a(\mu)\ge a_0+c_*\mu\qquad\text{for every }\mu>0,
		\end{equation}
		and $\cz\le c_*\le\rho_0$.
	\end{lemma}
	\begin{proof}
		Let $0<\mu_1<\mu_2$. By Proposition \ref{prop:fixed}\ref{it:compact},
		the function $c_+$ is upper semicontinuous, and hence it is bounded on
		$[\mu_1,\mu_2]$ by some $L>0$. By \eqref{eq:derivatives}, the
		continuous function $\mu\mapsto L\mu-a(\mu)$ has the left derivative
		$L-c_+(\mu)\ge0$ at every $\mu\in[\mu_1,\mu_2]$. A continuous function
		whose left derivative is nonnegative at every point is nondecreasing,
		see the proof of \cite[Lemma~3.4]{MZ12}. Hence $L\mu-a(\mu)$ is
		nondecreasing on $[\mu_1,\mu_2]$. Since $a$ is nondecreasing as well,
		we get
		\begin{equation*}
			0\le a(\mu)-a(\nu)\le L(\mu-\nu)
			\qquad\text{for }\mu_1\le\nu<\mu\le\mu_2 .
		\end{equation*}
		Thus $a$ is locally Lipschitz continuous on $(0,\infty)$, hence locally
		absolutely continuous and differentiable almost everywhere, and
		$a'=c_-$ at every point of differentiability by
		\eqref{eq:derivatives}. Since $c_-$ is lower semicontinuous, hence
		measurable, we get
		\begin{equation*}
			a(\mu)-a(\nu)=\int_\nu^\mu c_-(t)\,dt\qquad\text{for }0<\nu<\mu .
		\end{equation*}
		Since $a(\nu)\to a_0$ as $\nu\to0^+$, letting $\nu\to0^+$ and using the
		monotone convergence theorem we obtain \eqref{eq:integral0}. In
		particular, the integral is finite.
		Thus $\rho(\mu)$ is the mean value of $c_-$ over $(0,\mu)$. Since
		$c_-\ge c_*$ on $(0,\infty)$, we get $\rho(\mu)\ge c_*$ for every
		$\mu>0$, which is \eqref{eq:gap}, and hence $\rho_0\ge c_*$. Finally,
		$\cz\le c_*$ because $\cC(\mu)\subset\cM_\mu$ for every $\mu>0$.
	\end{proof}
	
	\begin{proof}[Proof of Theorem \ref{thm:action}]
		We first prove \eqref{eq:mainthresholds}. By Lemmas
		\ref{prop:zero-level} and \ref{prop:positive}, the limit $a_0$ exists,
		$0<a_0\le a(\mu)<\infty$ for every $\mu>0$, and $\cz>0$. By Lemma
		\ref{prop:thresholds}, $\cz\le c_*$, and $c_*<\infty$ because
		$\cC(\mu)\ne\varnothing$. The strict inequality $\cz<c_*$ is proved in
		Lemma \ref{prop:strict-smallmass} below.
		
		\ref{it:none} If $(\mu,u)$ is a solution of \eqref{eq:normalized}, then
		$u\in\cM_\mu$ by the Poho\v{z}aev identity, and
		\eqref{eq:smallmass-strict} gives $m^2=\frac12\norm[2]{u}^2>\cz$.
		
		\ref{it:below} If $u\in\cC(\mu)$, then
		$\frac12\norm[2]{u}^2\ge c_-(\mu)\ge c_*$.
		
		\ref{it:above} If $m^2>c_*$, then $c_-(\mu_1)<m^2$ for some $\mu_1>0$.
		By \eqref{eq:b} and \eqref{eq:derivatives}, $(b_m)'_+(\mu_1)=c_-(\mu_1)-m^2<0$, so that
		$b_m(\mu)<b_m(\mu_1)$ for $\mu>\mu_1$ close to $\mu_1$. Since $b_m$ is
		continuous and $b_m(\mu)=\mu\big(a(\mu)/\mu-m^2\big)\to+\infty$ as
		$\mu\to\infty$ by Proposition \ref{prop:fixed}\ref{it:large}, the
		minimum of $b_m$ on $[\mu_1,\infty)$ is attained at some $\mu>\mu_1$.
		Thus $\mu$ is a local minimum point of $b_m$, and the assertion follows
		from Proposition \ref{prop:fixed}\ref{it:selection}.
		
		\ref{it:at} If $c_-(\mu_0)=c_*$, then Proposition
		\ref{prop:fixed}\ref{it:compact} gives $u\in\cC(\mu_0)$ with
		$\frac12\norm[2]{u}^2=c_*$, and
		$\cE(u)=I_{\mu_0}(u)-\mu_0c_*=a(\mu_0)-\mu_0c_*\ge a_0$ by
		\eqref{eq:gap}. Conversely, if $u\in\cC(\mu_0)$ and
		$\frac12\norm[2]{u}^2=c_*$, then
		$c_*\le c_-(\mu_0)\le\frac12\norm[2]{u}^2=c_*$.
		
		The form of the elements of $\cC(\mu)$ is given by Proposition
		\ref{prop:fixed}\ref{it:exist}. If $f$ is odd, then $I_\mu$ and $M_\mu$
		are even, and $-u\in\cC(\mu)$ whenever $u\in\cC(\mu)$.
	\end{proof}
	
	\section{Energy ground states and least energy normalized solutions}\label{sec:energy}

	We first prove Theorem \ref{thm:energy}, using neither
	\Fthree\ nor Proposition \ref{prop:fixed}. Since $\Psi$ is not
	differentiable, the Lagrange multiplier of a constrained local minimum
	point is obtained from Theorem \ref{thm:multiplier} in Appendix
	\ref{app:multiplier}.
	
	\begin{lemma}\label{lem:em-bounds}
		Assume \Fone--\Ftwo. Then the levels $e_m$ defined in \eqref{eq:mass}
		satisfy
		\begin{equation}\label{eq:em-bounds}
			-\infty<e_m\le0\qquad\text{for every }m>0 .
		\end{equation}
		Moreover, for every $u\in\D\setminus\{0\}$ and $\theta>1$,
		\begin{equation}\label{eq:direct-strict-scaling}
			\cE(u_\theta)=\theta^N\Psi(u)-\theta^N\int_{\RN}F(\theta u)\dd
			\le\theta^{N+2}\cE(u)-(\theta^{N+2}-\theta^N)\Psi(u)
			<\theta^{N+2}\cE(u).
		\end{equation}
	\end{lemma}
	\begin{proof}
		We have $e_m>-\infty$ by \cite[Lemma~A.2]{MZ12}. Let
		$\phi\in\cC_c^\infty(\RN)$ with $\norm[2]{\phi}^2=2m^2$. For large $R$
		the function $v=R^{-N/2}\phi(\cdot/R)$ belongs to $\cS_m$, and by Lemma
		\ref{lem:f}\ref{it:f-sign} and \eqref{eq:kinetic},
		\begin{equation*}
			e_m\le\cE(v)\le\Psi(v)\le\norm[2]{\nabla v}^2
			=R^{-2}\norm[2]{\nabla\phi}^2 ,
		\end{equation*}
		which tends to zero as $R\to\infty$. This proves \eqref{eq:em-bounds}.
		In \eqref{eq:direct-strict-scaling}, the equality follows from
		\eqref{eq:scaling}, the first inequality from Lemma
		\ref{lem:f}\ref{it:f-mono}, and the second one from $\Psi(u)>0$.
	\end{proof}
	
	\begin{lemma}\label{lem:direct-multiplier}
		Assume \Fone--\Ftwo, let $m>0$, and let $u\in\cS_m$ be a local minimum
		point of $\cE$ on $\cS_m$ in the topology of $X$. Then $u$ satisfies
		the variational inequality \eqref{eq:critical} with
		\begin{equation}\label{eq:direct-multiplier}
			\mu=\mu_u:=\frac{N}{(N+2)m^2}
			\left(\int_{\RN}G(u)\dd-\Psi(u)\right).
		\end{equation}
		If, in addition, $\cE(u)\le0$, then $\mu_u>0$ and $(\mu_u,u)$ is a
		solution of \eqref{eq:normalized}. In this case $u$ is a classical
		solution of the equation in \eqref{eq:normalized} with $\mu=\mu_u$, and
		$\norm[\infty]{\nabla u}<1$.
	\end{lemma}
	\begin{proof}
		We apply Theorem \ref{thm:multiplier} to the functionals $\Psi$,
		$\Phi(w)=\int_{\RN}F(w)\dd$ and $Q(w)=\frac12\norm[2]{w}^2$ on $X$,
		with $q=m^2$ and with the dilations
		$T_\theta w=\theta w(\cdot/\theta)$, $\frac12<\theta<2$. The operators
		$T_\theta$ are linear, and
		\begin{equation*}
			\norm[2]{T_\theta w}=\theta^{(N+2)/2}\norm[2]{w},\qquad
			\norm[r]{\nabla T_\theta w}=\theta^{N/r}\norm[r]{\nabla w}
			\quad(r=2,2j_0).
		\end{equation*}
		Hence they are uniformly bounded on $X$, and $T_\theta w\to w$ in $X$
		as $\theta\to1$ by the continuity of dilations in Lebesgue spaces.
		By \eqref{eq:scaling}, $\Psi(T_\theta w)=\theta^N\Psi(w)$,
		$Q(T_\theta w)=\theta^{N+2}Q(w)$ and
		$\Phi(T_\theta w)=\theta^N\Phi(\theta w)$. Since $\Phi$ and $Q$ are of
		class $\cC^1$ on $X$, the assumptions of Theorem \ref{thm:multiplier}
		hold with $\gamma(\theta)=\theta^N$, and
		\begin{equation*}
			\frac{d}{d\theta}Q(T_\theta u)\Big|_{\theta=1}=(N+2)m^2,\qquad
			\frac{d}{d\theta}\Phi(T_\theta u)\Big|_{\theta=1}
			=N\int_{\RN}F(u)\dd+\int_{\RN}f(u)u\dd=N\int_{\RN}G(u)\dd .
		\end{equation*}
		Thus the number $\lambda$ in \eqref{eq:multiplier} equals $\mu_u$, and
		\eqref{eq:multiplier-rule} is \eqref{eq:critical} with $\mu=\mu_u$.
		
		If $\cE(u)\le0$, then
		$\int_{\RN}G(u)\dd-\Psi(u)=-\cE(u)+\frac1N\int_{\RN}f(u)u\dd>0$ by
		Lemma \ref{lem:f}\ref{it:f-sign}, hence $\mu_u>0$ and $u$ is a critical
		point of $I_{\mu_u}$. In this case the function
		$\sqrt{\mu_u}\,u(\cdot/\sqrt{\mu_u})$ is a critical point of the
		functional associated with the frequency one and with the nonlinearity
		$g(s)=\mu_u^{-1/2}f(s/\sqrt{\mu_u})$, which is of class $\cC^1$ with
		$g(0)=g'(0)=0$. Hence the positive mass part of
		\cite[Lemma~2.3]{BIM}
		shows that $u$ is a classical solution with
		$\norm[\infty]{\nabla u}<1$.
	\end{proof}
	
	\begin{proof}[Proof of Theorem \ref{thm:energy}]
		\emph{Step 1: the threshold.} We claim that
		$0<\widetilde r_e<\infty$ and
		\begin{equation}\label{eq:em-sign}
			e_m=0\quad\text{for }m^2\le\widetilde r_e,\qquad
			e_m<0\quad\text{for }m^2>\widetilde r_e .
		\end{equation}
		Let $u\in\cS_m$ with $m\le1$. Then $\norm[\infty]{u}\le C$ by Lemma
		\ref{lem:amplitude}, and Lemma \ref{lem:f}\ref{it:f-local} with
		$q=2+4/N$, the Gagliardo--Nirenberg inequality and \eqref{eq:kinetic}
		give
		\begin{equation*}
			\cE(u)\ge\frac12\norm[2]{\nabla u}^2-C\norm[2+4/N]{u}^{2+4/N}
			\ge\Big(\frac12-Cm^{4/N}\Big)\norm[2]{\nabla u}^2 .
		\end{equation*}
		Hence $e_m\ge0$ for small $m$, so that $\widetilde r_e>0$. On the other
		hand, if $\phi\in\cC_c^\infty(\RN)\setminus\{0\}$ and
		$\norm[\infty]{\nabla\phi}\le1$, then \eqref{eq:direct-strict-scaling}
		gives
		$\cE(\phi_\theta)\le\theta^N\Psi(\phi)-\theta^{N+2}\int_{\RN}F(\phi)\dd<0$
		for large $\theta$, hence $\widetilde r_e<\infty$. Let now $e_m<0$, and
		choose $u\in\cS_m$ with $\cE(u)<0$. By
		\eqref{eq:direct-strict-scaling}, $\cE(u_\theta)<0$ for $\theta>1$, and
		by the continuity of $\theta\mapsto\cE(u_\theta)$ also for $\theta<1$
		close to $1$. Since
		$\frac12\norm[2]{u_\theta}^2=\theta^{N+2}m^2$, the set
		$\{m^2>0:e_m<0\}$ is open and, together with each of its elements, it
		contains all larger numbers. Hence this set is
		$(\widetilde r_e,\infty)$, and \eqref{eq:em-sign} follows from
		\eqref{eq:em-bounds}.
		
		\emph{Step 2: reduction to radial functions which do not change sign.}
		Let $u\in\cS_m$. We claim that there is $v\in\cS_m$ of the form
		\begin{equation}\label{eq:direct-radial}
			v(x)=\pm V(|x|),\qquad\text{where }V\ge0\text{ is nonincreasing on }
			[0,\infty),
		\end{equation}
		such that $\cE(v)\le\cE(u)$, with strict inequality if $u$ changes
		sign. Suppose first that $u$ changes sign, and let $u^+=\max\{u,0\}$
		and $u^-=\min\{u,0\}$. Then $u^\pm\in\D\setminus\{0\}$,
		$\cE(u)=\cE(u^+)+\cE(u^-)$ and
		$\norm[2]{u}^2=\norm[2]{u^+}^2+\norm[2]{u^-}^2$. Hence one of the
		functions $u^+$ and $u^-$, which we denote by $z$, satisfies
		$\cE(z)/\norm[2]{z}^2\le\cE(u)/\norm[2]{u}^2$. Let $\theta>1$ be such
		that $\theta^{N+2}\norm[2]{z}^2=\norm[2]{u}^2$. Then $z_\theta\in\cS_m$,
		and \eqref{eq:direct-strict-scaling} gives
		\begin{equation*}
			\cE(z_\theta)<\theta^{N+2}\cE(z)\le\cE(u).
		\end{equation*}
		Therefore we may assume that $u$ does not change sign, and, replacing
		$f$ by $s\mapsto-f(-s)$ and $u$ by $-u$ if necessary, that $u\ge0$.
		Let $u^*$ be the Schwarz symmetrization of $u$. Since $u^*$ and $u$ are
		equimeasurable, we have $\norm[2]{u^*}=\norm[2]{u}$ and
		$\int_{\RN}F(u^*)\dd=\int_{\RN}F(u)\dd$. Moreover,
		$u\in W^{1,2j}(\RN)$ for every integer $j\ge1$, because $u$ and
		$\nabla u$ belong to $L^2(\RN)\cap L^\infty(\RN)$, and the
		P\'olya--Szeg\H{o} inequality gives
		\begin{equation*}
			\norm[2j]{\nabla u^*}\le\norm[2j]{\nabla u}
			\le\norm[2]{\nabla u}^{1/j} .
		\end{equation*}
		Letting $j\to\infty$ we get $\norm[\infty]{\nabla u^*}\le1$, hence
		$u^*\in\cS_m$. Finally, since
		$1-\sqrt{1-t^2}=\sum_{j\ge1}b_jt^{2j}$ for $0\le t\le1$ with $b_j>0$,
		the same inequalities and the monotone convergence theorem give
		$\Psi(u^*)\le\Psi(u)$, see also the proof of \cite[Lemma~B.2]{MZ12}.
		Thus $\cE(u^*)\le\cE(u)$. Since $u^*$ is of the form
		\eqref{eq:direct-radial}, the claim follows.
		
		\emph{Step 3: compactness.} Let $(u_n)\subset\D$ be a sequence of
		functions of the form \eqref{eq:direct-radial} such that
		$\norm[2]{u_n}$ and $\Psi(u_n)$ are bounded. We claim that, up to a
		subsequence, $u_n\to u$ locally uniformly, where $u\in\D$ is of the
		form \eqref{eq:direct-radial}, and
		\begin{equation}\label{eq:direct-compactness}
			\norm[2]{u}\le\liminf_{n\to\infty}\norm[2]{u_n},\qquad
			\Psi(u)\le\liminf_{n\to\infty}\Psi(u_n),\qquad
			\int_{\RN}F(u_n)\dd\to\int_{\RN}F(u)\dd .
		\end{equation}
		Indeed, we may assume that the functions $u_n$ have the same sign. By
		Lemma \ref{lem:amplitude} they are uniformly bounded, and they are
		$1$-Lipschitz. By the Arzel\`a--Ascoli theorem, up to a subsequence,
		$u_n\to u$ locally uniformly, where $u$ is $1$-Lipschitz and of the
		form \eqref{eq:direct-radial}, and the first inequality in
		\eqref{eq:direct-compactness} follows from Fatou's lemma. By
		\eqref{eq:kinetic}, $(\nabla u_n)$ is bounded in $L^2(\RN)$, hence, up
		to a subsequence, it converges weakly in $L^2(\RN)$, and its limit is
		$\nabla u$ because $u_n\to u$ locally uniformly. In particular,
		$u\in H^1(\RN)$ and $|\nabla u|\le1$, so that $u\in\D$. The functional
		$g\mapsto\int_{\RN}\big(1-\sqrt{1-|g|^2}\big)\dd$ is convex and, by
		Fatou's lemma, lower semicontinuous on the closed convex subset
		$\{g:|g|\le1\}$ of $L^2(\RN;\RN)$. Hence it is weakly lower
		semicontinuous there, which gives the second inequality in
		\eqref{eq:direct-compactness}. Finally, the radial monotonicity gives
		$|u_n(x)|\le C\min\{1,|x|^{-N/2}\}$, and Lemma
		\ref{lem:f}\ref{it:f-local} with $q=2^*$ gives
		$F(u_n)\le C\min\{1,|x|^{-N^2/(N-2)}\}$. Thus the last assertion in
		\eqref{eq:direct-compactness} follows from the dominated convergence
		theorem.
		
		\emph{Step 4: the case $m^2>\widetilde r_e$.} By \eqref{eq:em-sign},
		$e_m<0$. By Step 2 there is a minimizing sequence $(u_n)\subset\cS_m$
		of functions of the form \eqref{eq:direct-radial}. Since the functions
		in $\cS_m$ are uniformly bounded, Lemma \ref{lem:f}\ref{it:f-local} shows that
		$\Psi(u_n)=\cE(u_n)+\int_{\RN}F(u_n)\dd$ is bounded. Let $u$ be given
		by Step 3. Then $\frac12\norm[2]{u}^2\le m^2$ and
		$\cE(u)\le\liminf_{n\to\infty}\cE(u_n)=e_m<0$. In particular, $u\ne0$.
		If $\frac12\norm[2]{u}^2<m^2$, then $u_\theta\in\cS_m$ for some
		$\theta>1$, and \eqref{eq:direct-strict-scaling} gives
		\begin{equation*}
			e_m\le\cE(u_\theta)<\theta^{N+2}\cE(u)<\cE(u)\le e_m,
		\end{equation*}
		a contradiction. Hence $u\in\cS_m$ and $\cE(u)=e_m$, which proves
		\ref{it:en-above}.
		
		\emph{Step 5: the case $m^2\le\widetilde r_e$.} Let
		$v\in\D\setminus\{0\}$ with $\cE(v)\le0$. By
		\eqref{eq:direct-strict-scaling}, $\cE(v_\theta)<0$ for every
		$\theta>1$, hence $\theta^{N+2}\frac12\norm[2]{v}^2>\widetilde r_e$ by
		\eqref{eq:em-sign}. Letting $\theta\to1^+$ we get
		\begin{equation}\label{eq:least-mass}
			\tfrac12\norm[2]{v}^2\ge\widetilde r_e
			\qquad\text{for every }v\in\D\setminus\{0\}\text{ with }\cE(v)\le0 .
		\end{equation}
		In particular, if $m^2<\widetilde r_e$, then $\cE>0=e_m$ on $\cS_m$,
		which is \ref{it:en-below}. Next, by \eqref{eq:em-sign} there are
		$u_n\in\D\setminus\{0\}$ such that $\cE(u_n)<0$ and
		$\frac12\norm[2]{u_n}^2\to\widetilde r_e$, and by Step 2 we may assume
		that they are of the form \eqref{eq:direct-radial}. Then $(u_n)$ is
		bounded in $L^\infty(\RN)$ by Lemma \ref{lem:amplitude}, and Lemma
		\ref{lem:f}\ref{it:f-local} with $q=2$ and $q=2^*$, the Sobolev
		inequality and \eqref{eq:kinetic} give
		\begin{equation*}
			\Psi(u_n)<\int_{\RN}F(u_n)\dd
			\le C\min\big\{\norm[2]{u_n}^2,\norm[2^*]{u_n}^{2^*}\big\}
			\le C\min\big\{1,\Psi(u_n)^{2^*/2}\big\}.
		\end{equation*}
		Since $2^*>2$, this shows that $\Psi(u_n)$, and hence also
		$\int_{\RN}F(u_n)\dd$, is bounded and bounded away from zero. Let
		$u$ be given by Step 3. Then $\int_{\RN}F(u)\dd>0$, so that $u\ne0$,
		and
		\begin{equation*}
			\tfrac12\norm[2]{u}^2\le\widetilde r_e,\qquad
			\cE(u)\le\liminf_{n\to\infty}\cE(u_n)\le0 .
		\end{equation*}
		By \eqref{eq:least-mass}, $\frac12\norm[2]{u}^2=\widetilde r_e$. Thus
		the minimum in \eqref{eq:threshold-minimum-mass} is attained and equals
		$\widetilde r_e$. Moreover, for $m^2=\widetilde r_e$ we have
		$u\in\cS_m$ and $\cE(u)\le0=e_m$, which proves \ref{it:en-at}.
		
		\emph{Step 6: properties of the energy ground states.} Let $u\in\cS_m$
		satisfy $\cE(u)=e_m$. Then $\cE(u)\le0$, and $u$ does not change sign
		by Step 2. By Lemma \ref{lem:direct-multiplier}, $(\mu_u,u)$ is a
		solution of \eqref{eq:normalized} with $\mu_u>0$, and $u$ is of class
		$\cC^2$ with $\norm[\infty]{\nabla u}<1$. Hence the equation in
		\eqref{eq:normalized} can be written as
		\begin{equation*}
			\sum_{i,j=1}^Na_{ij}\,\partial_i\partial_ju-\mu_uu=-f(u),\qquad
			a_{ij}:=\frac{\delta_{ij}}{\sqrt{1-|\nabla u|^2}}
			+\frac{\partial_iu\,\partial_ju}{(1-|\nabla u|^2)^{3/2}},
		\end{equation*}
		where the coefficients $a_{ij}$ are continuous, bounded and uniformly
		elliptic. If $u\ge0$, then $f(u)\ge0$ by Lemma
		\ref{lem:f}\ref{it:f-sign}, and the strong maximum principle
		\cite[Theorem~3.5]{GT} gives $u>0$ in $\RN$. Similarly, $u<0$ in $\RN$
		if $u\le0$. Finally, the energy ground states obtained in Steps 4 and 5
		are of the form \eqref{eq:direct-radial}, and $V>0$ by the above.
	\end{proof}

	The following proposition does not depend on \Ffour. For $2<p<2^*$ one
	has $B_m\le0$, and it reduces to
	\cite[Theorem~1.3]{MZ12}, whose proof we follow.
	
	\begin{proposition}\label{prop:duality}
		Assume \Fone--\Fthree. Then
		\begin{equation}\label{eq:duality}
			e_m=\min\{0,B_m\}\qquad\text{for every }m>0,
		\end{equation}
		and the function $m^2\mapsto e_m$ is concave on $(0,\infty)$. Moreover,
		for every $m>0$
		the level $e_m$ is attained if and only if $B_m\le0$ and $b_m$ attains
		its infimum on $(0,\infty)$. In this case every global minimum point
		$\mu$ of $b_m$ satisfies $\cC(\mu)\subset\cS_m$, and
		\begin{equation}\label{eq:minimizers}
			\{u\in\cS_m:\cE(u)=e_m\}
			=\bigcup_{\substack{\mu>0\\ b_m(\mu)=B_m}}\cC(\mu).
		\end{equation}
	\end{proposition}
	\begin{proof}
		Let $m>0$, $\mu>0$, $u\in\cC(\mu)$ and
		$\theta=(2m^2/\norm[2]{u}^2)^{1/(N+2)}$.
		Then $u_\theta\in\cS_m$ by \eqref{eq:scaling}, and Proposition
		\ref{prop:fixed}\ref{it:fibres} gives
		\begin{equation}\label{eq:firstcomparison}
			e_m\le\cE(u_\theta)=I_\mu(u_\theta)-\mu m^2
			\le I_\mu(u)-\mu m^2=b_m(\mu),
		\end{equation}
		where the second inequality is strict unless $\theta=1$.
		Together with \eqref{eq:em-bounds}, this shows that
		$e_m\le\min\{0,B_m\}$.
		
		On the other hand, let $u\in\cS_m$ satisfy $\cE(u)\le0$, and let
		$\mu_u$ be given by \eqref{eq:direct-multiplier}. Since
		$\int_{\RN}F(u)\dd\ge\Psi(u)$ and $\int_{\RN}f(u)u\dd>0$, we have
		$\mu_u>0$, and \eqref{eq:M} gives $u\in\cM_{\mu_u}$.
		Hence
		\begin{equation}\label{eq:secondcomparison}
			\cE(u)=I_{\mu_u}(u)-\mu_um^2\ge a(\mu_u)-\mu_um^2=b_m(\mu_u)\ge B_m,
		\end{equation}
		where the first inequality is an equality if and only if
		$u\in\cC(\mu_u)$. Consequently, $\cE(u)\ge\min\{0,B_m\}$ for every
		$u\in\cS_m$, which proves \eqref{eq:duality}. The concavity
		follows, because $e_m$ is the infimum of the constant zero and of the
		affine functions $m^2\mapsto a(\mu)-\mu m^2$.
		
		If $u\in\cS_m$ and $\cE(u)=e_m$, then $\cE(u)\le0$ and, by
		\eqref{eq:secondcomparison},
		\begin{equation*}
			e_m=\cE(u)\ge b_m(\mu_u)\ge B_m\ge e_m.
		\end{equation*}
		Thus $B_m=e_m\le0$, the frequency $\mu_u$ minimizes $b_m$, and
		$u\in\cC(\mu_u)$.
		Conversely, suppose that $B_m\le0$ and $b_m(\mu)=B_m$, so that
		$e_m=B_m$. For every $u\in\cC(\mu)$ both inequalities in
		\eqref{eq:firstcomparison} are equalities, and hence $\theta=1$.
		Thus $u\in\cS_m$ and $\cE(u)=e_m$.
		This proves the criterion and \eqref{eq:minimizers}.
	\end{proof}
	
	Under \Fone--\Ffour, Lemmas \ref{prop:zero-level} and
	\ref{prop:positive} and Proposition \ref{prop:fixed}\ref{it:large}
	give
	\begin{equation}\label{eq:blimits}
		\lim_{\mu\to0^+}b_m(\mu)=a_0>0,\qquad
		\lim_{\mu\to\infty}b_m(\mu)=+\infty.
	\end{equation}
	
	\begin{proof}[Proof of Theorem \ref{thm:least}]
		By Lemma \ref{prop:thresholds}, $c_*\le\rho_0$. To prove that
		$\rho_0<\widetilde r_e$, let $m^2=\widetilde r_e$. By Theorem
		\ref{thm:energy}\ref{it:en-at}, the level $e_m=0$ is attained. Hence,
		by Proposition \ref{prop:duality}, $B_m\le0$ and $b_m$ attains its
		infimum, whereas $\min\{0,B_m\}=e_m=0$. Thus $b_m(\mu)=B_m=0$ for some
		$\mu>0$, that is, $a(\mu)=\mu\widetilde r_e$, and
		\begin{equation*}
			\rho_0\le\rho(\mu)=\widetilde r_e-\frac{a_0}{\mu}<\widetilde r_e ,
		\end{equation*}
		because $a_0>0$.
		
		\ref{it:least-lower} Let $(\mu,u)$ be a solution of
		\eqref{eq:normalized}. Since $u\in\cM_\mu$ by the Poho\v{z}aev
		identity, we have $I_\mu(u)\ge a(\mu)$ and
		\begin{equation}\label{eq:solutionsbound}
			\cE(u)=I_\mu(u)-\mu m^2\ge b_m(\mu)\ge B_m,
		\end{equation}
		where the first inequality is an equality if and only if
		$u\in\cC(\mu)$. If $m^2<\rho_0$, then
		$b_m(\mu)\ge a_0+\mu(\rho_0-m^2)>a_0$ by \eqref{eq:q}. Let now
		$m^2<\widetilde r_e$, and suppose that $B_m\le0$. By \eqref{eq:blimits}
		and the continuity of $b_m$, the infimum $B_m$ is attained. Hence $e_m$
		is attained by Proposition \ref{prop:duality}, which contradicts
		Theorem \ref{thm:energy}\ref{it:en-below}. Thus $B_m>0$.
		
		\ref{it:least-attain} By \eqref{eq:blimits}, $B_m\le a_0$. If
		$m^2>\rho_0$, then $\rho(\mu)<m^2$ for some $\mu>0$, hence
		$B_m\le b_m(\mu)<a_0$ by \eqref{eq:q}, and by \eqref{eq:blimits} and
		the continuity of $b_m$ the infimum of $b_m$ is attained. If
		$m^2<\rho_0$, then $b_m>a_0\ge B_m$ on
		$(0,\infty)$, and the infimum is not attained. If $m^2=\rho_0$, then
		$b_m\ge a_0$, so that $B_m=a_0$, and $b_m(\mu)=B_m$ if and only if
		$\rho(\mu)=\rho_0$. This proves the first statement.
		Suppose now that $b_m(\mu_m)=B_m$. By Proposition
		\ref{prop:fixed}\ref{it:selection}, $\cC(\mu_m)\subset\cS_m$. Hence, for
		every $u\in\cC(\mu_m)$, the pair $(\mu_m,u)$ is a solution of
		\eqref{eq:normalized} with $\cE(u)=b_m(\mu_m)=B_m$, and by
		\eqref{eq:solutionsbound} it has the least energy among the solutions
		of \eqref{eq:normalized}. Conversely, if $(\mu,u)$ is a solution with
		$\cE(u)=B_m$, then equality holds throughout
		\eqref{eq:solutionsbound}, that is, $u\in\cC(\mu)$ and $b_m(\mu)=B_m$.
	\end{proof}
	
	\section{The zero-frequency limit and attainment of the mass thresholds}\label{sec:endpoints}
	
	Throughout this section we assume \Fone--\Ffour. Replacing $f$ by
	$s\mapsto-f(-s)$ and $u$ by $-u$ preserves \Fone--\Ffour, the relations
	$u\in\cC(\mu)$ and $u\in\cM_0$, and the quantities $a(\mu)$, $a_0$, $\cz$
	and $c_\pm(\mu)$. Therefore, by Proposition
	\ref{prop:fixed}\ref{it:exist}, given $\mu_n>0$ and $u_n\in\cC(\mu_n)$,
	after translations, passing to a subsequence, and this replacement if
	necessary, we may assume that
	\begin{equation}\label{eq:positive-radial}
		u_n(x)=W_n(|x|),\qquad\text{where $W_n$ is positive and decreasing on }
		[0,\infty).
	\end{equation}
	
	Let $u(x)=W(|x|)\in\cC(\mu)$ be positive. Since $u$ is a classical
	solution, integrating the equation in \eqref{eq:normalized} over the
	ball $\{|x|<r\}$ we get
	\begin{equation}\label{eq:flux}
		W'(r)=-\frac{Z(r)}{\sqrt{1+Z(r)^2}},\qquad
		Z(r):=r^{1-N}\int_0^r s^{N-1}\big(f(W(s))-\mu W(s)\big)\,ds
		\qquad(r>0).
	\end{equation}
	Since $W'\le0$, we have $Z\ge0$, and
	\begin{equation}\label{eq:fluxbound}
		0\le-W'(r)\le Z(r)\le r^{1-N}\int_0^rs^{N-1}f(W(s))\,ds\qquad(r>0).
	\end{equation}
	We also recall that, by \eqref{eq:functionals} and \eqref{eq:M},
	\begin{equation}\label{eq:actionidentity}
		a(\mu)=I_\mu(u)=\frac1N\left(\int_{\RN}f(u)u\dd-\mu\norm[2]{u}^2\right)
		\qquad\text{for }u\in\cC(\mu).
	\end{equation}
	
	\begin{lemma}\label{lem:locallimit}
		Let $\mu_n\to0^+$ and $u_n\in\cC(\mu_n)$ satisfy
		\eqref{eq:positive-radial}, and suppose that $(u_n)$ is bounded in
		$L^2(\RN)$. Then, up to a subsequence, $W_n\to W$ in
		$\cC^1([0,R])$ for every $R>0$, where $W\in\cC^1([0,\infty))$ has the
		following properties.
		\begin{enumerate}[label=\textnormal{(\alph*)}]
			\item\label{it:ll-basic}
			$W$ is positive and nonincreasing, $W(r)\to0$ as $r\to\infty$, and the
			function $w(x):=W(|x|)$ satisfies
			$\norm[2]{w}^2\le\liminf_{n\to\infty}\norm[2]{u_n}^2$.
			\item\label{it:ll-eq}
			$w$ is a classical solution of \eqref{eq:zerofreq}, and
			$\norm[\infty]{\nabla w}<1$.
			\item\label{it:ll-lower}
			There is $c>0$ such that $W(r)\ge c\,r^{2-N}$ for $r\ge1$.
		\end{enumerate}
	\end{lemma}
	\begin{proof}
		By Lemma \ref{lem:amplitude}, the sequence $(u_n)$ is bounded in
		$L^\infty(\RN)$, and the functions $W_n$ are $1$-Lipschitz. By the
		Arzel\`a--Ascoli theorem and a diagonal argument, up to a subsequence,
		$W_n\to W$ uniformly on compact subsets of $[0,\infty)$, where $W$ is
		nonnegative, nonincreasing and $1$-Lipschitz. We have $W(0)>0$.
		Otherwise $\norm[\infty]{u_n}=W_n(0)\to0$, and
		\eqref{eq:actionidentity}, Lemma
		\ref{prop:zero-level}\ref{it:zl-limit}, Lemma \ref{prop:positive}
		and \Fone\ would give
		\begin{equation*}
			0<a_0\le a(\mu_n)\le\frac1N\int_{\RN}f(u_n)u_n\dd
			\le\frac{\norm[2]{u_n}^2}{N}\sup_{0<s\le W_n(0)}\frac{f(s)}s\to0 .
		\end{equation*}
		
		Let $Z_n$ be the function associated with $W_n$ and $\mu_n$ by
		\eqref{eq:flux}, and set
		\begin{equation}\label{eq:zeroflux}
			Z(r):=r^{1-N}\int_0^r s^{N-1}f(W(s))\,ds\quad(r>0),\qquad Z(0):=0 .
		\end{equation}
		Since $r^{1-N}\int_0^rs^{N-1}\,ds=r/N$, for every $R>0$ we have
		\begin{equation*}
			\sup_{0<r\le R}|Z_n(r)-Z(r)|
			\le\frac RN\sup_{0\le s\le R}
			\big(|f(W_n(s))-f(W(s))|+\mu_nW_n(s)\big)\to0 .
		\end{equation*}
		By \eqref{eq:flux}, $W_n'$ converges uniformly on $[0,R]$ to
		$-Z/\sqrt{1+Z^2}$. Hence $W\in\cC^1([0,\infty))$, $W_n\to W$ in
		$\cC^1([0,R])$, and
		\begin{equation}\label{eq:zeroprofile}
			W'=-\frac{Z}{\sqrt{1+Z^2}}\qquad\text{on }[0,\infty).
		\end{equation}
		Since $Z'=f(W)-(N-1)Z/r$ on $(0,\infty)$ and $Z(r)/r\to f(W(0))/N$ as
		$r\to0^+$, we have $W\in\cC^2((0,\infty))$ and
		\begin{equation*}
			\lim_{r\to0^+}\frac{W'(r)}r
			=\lim_{r\to0^+}W''(r)=-\frac{f(W(0))}{N}.
		\end{equation*}
		Thus $w\in\cC^2(\RN)$. Moreover, \eqref{eq:zeroflux} and
		\eqref{eq:zeroprofile} give $Z=-W'/\sqrt{1-W'^2}$ and
		$(r^{N-1}Z)'=r^{N-1}f(W)$, that is, $w$ is a classical solution of
		\eqref{eq:zerofreq}.
		
		Fatou's lemma gives
		$\norm[2]{w}^2\le\liminf_{n\to\infty}\norm[2]{u_n}^2$. Since $W$ is
		nonincreasing, we have $\frac{\omega_N}{N}r^NW(r)^2\le\norm[2]{w}^2$
		for $r>0$, and in particular $W(r)\to0$ as $r\to\infty$. Moreover, by
		Lemma \ref{lem:f}\ref{it:f-local} with $q=2^*$,
		$f(W(s))\le CW(s)^{2^*-1}\le Cs^{-N(N+2)/(2(N-2))}$ for $s\ge1$, and
		the exponent is larger than one. Hence
		\begin{equation*}
			Z(r)\le\int_0^rf(W(s))\,ds\le\int_0^\infty f(W(s))\,ds<\infty ,
		\end{equation*}
		so that $Z$ is bounded and $\norm[\infty]{\nabla w}=\sup_{r\ge0}|W'(r)|<1$
		by \eqref{eq:zeroprofile}.
		
		Finally, since $W(0)>0$, the number
		$c:=\int_0^1s^{N-1}f(W(s))\,ds$ is positive, and $Z(r)\ge c\,r^{1-N}$
		for $r\ge1$. Since $t\mapsto t/\sqrt{1+t^2}$ is increasing,
		\eqref{eq:zeroprofile} gives
		\begin{equation*}
			-W'(r)\ge\frac{c\,r^{1-N}}{\sqrt{1+c^2r^{2-2N}}}
			\ge\frac{c}{\sqrt{1+c^2}}\,r^{1-N}\qquad(r\ge1) .
		\end{equation*}
		Integrating over $(r,\infty)$ and using $W(r)\to0$ as $r\to\infty$, we
		obtain \ref{it:ll-lower}. In particular, $W>0$ on $[0,\infty)$.
	\end{proof}
	
	\begin{lemma}\label{prop:strict-smallmass}
		Assume \Fone--\Ffour. Then the numbers defined in \eqref{eq:thresholds}
		satisfy $\cz<c_*$.
	\end{lemma}
	\begin{proof}
		By Lemma \ref{prop:thresholds}, $\cz\le c_*$. Suppose by
		contradiction that $\cz=c_*$, and choose $\mu_n>0$ and
		$u_n\in\cC(\mu_n)$ such that
		$\frac12\norm[2]{u_n}^2=c_-(\mu_n)\to c_*$, see Proposition
		\ref{prop:fixed}\ref{it:compact}. By Proposition
		\ref{prop:fixed}\ref{it:large}, the sequence $(\mu_n)$ is bounded, and
		we may assume that $\mu_n\to\mu_0\ge0$. If $\mu_0>0$, then Proposition
		\ref{prop:fixed}\ref{it:compact} gives $y_n\in\RN$ and
		$u\in\cC(\mu_0)$ such that, up to a subsequence,
		$u_n(\cdot+y_n)\to u$ in $X$. Then
		$\frac12\norm[2]{u}^2=c_*=\cz$, which contradicts
		\eqref{eq:smallmass-strict}. Hence $\mu_0=0$.
		
		As explained at the beginning of this section, we may assume
		\eqref{eq:positive-radial}. Let $W$ and $w$ be given by Lemma
		\ref{lem:locallimit}. Then $w$ is a positive classical solution of
		\eqref{eq:zerofreq} with $\norm[\infty]{\nabla w}<1$ and
		\begin{equation}\label{eq:w-mass}
			\frac12\norm[2]{w}^2\le\lim_{n\to\infty}\frac12\norm[2]{u_n}^2=\cz .
		\end{equation}
		Since $W_n$ is nonincreasing and $(u_n)$ is bounded in $L^2(\RN)$ and
		in $L^\infty(\RN)$, we have $W_n(r)\le C\min\{1,r^{-N/2}\}$, and Lemma
		\ref{lem:f}\ref{it:f-local} with $q=2^*$ gives
		$G(u_n)\le C\min\{1,|x|^{-N^2/(N-2)}\}$. By the dominated convergence
		theorem, $\int_{\RN}G(u_n)\dd\to\int_{\RN}G(w)\dd$. Since
		$\nabla u_n\to\nabla w$ locally uniformly, Fatou's lemma and
		$M_{\mu_n}(u_n)=0$ give
		\begin{equation*}
			\begin{aligned}
				\int_{\RN}\big(1-\sqrt{1-|\nabla w|^2}\big)\dd
				&\le\liminf_{n\to\infty}\Psi(u_n)\\
				&=\lim_{n\to\infty}\left(\int_{\RN}G(u_n)\dd
				-\frac{N+2}{2N}\mu_n\norm[2]{u_n}^2\right)
				=\int_{\RN}G(w)\dd .
			\end{aligned}
		\end{equation*}
		Since $1-\sqrt{1-|\xi|^2}\ge\frac12|\xi|^2$ for $|\xi|\le1$, we get
		$\nabla w\in L^2(\RN)$. Together with \eqref{eq:w-mass} and
		$|\nabla w|\le1$, this shows that $w\in\D$ and $M_0(w)\le0$. By
		Lemma \ref{prop:zero-level}\ref{it:zl-fibre}, we have
		$\theta_0(w)\le1$, and \eqref{eq:smallmass-characterization} and
		\eqref{eq:w-mass} give
		\begin{equation*}
			\cz\le\frac12\norm[2]{w_{\theta_0(w)}}^2
			=\theta_0(w)^{N+2}\,\frac12\norm[2]{w}^2\le\frac12\norm[2]{w}^2\le\cz .
		\end{equation*}
		Hence $\theta_0(w)=1$, that is, $w\in\cM_0$, and $w$ minimizes
		$\frac12\norm[2]{u}^2$ among $u\in\cM_0$.
		
		Next, let $\varphi,\phi\in\cC_c^\infty(\RN)$ with $\phi\ge0$ and
		$\phi\ne0$. Since $\norm[\infty]{\nabla w}<1$, we have
		$w+t\varphi+\sigma\phi\in\D\setminus\{0\}$ for $|t|+|\sigma|$ small,
		and the function $(t,\sigma)\mapsto M_0(w+t\varphi+\sigma\phi)$ is of
		class $\cC^1$ near the origin. Testing
		\eqref{eq:zerofreq} with $\varphi$, we see that its partial derivative
		with respect to $t$ at the origin is
		\begin{equation*}
			\int_{\RN}\frac{\nabla w\cdot\nabla\varphi}{\sqrt{1-|\nabla w|^2}}\dd
			-\int_{\RN}G'(w)\varphi\dd
			=-\frac1N\int_{\RN}\big(wf'(w)+f(w)\big)\varphi\dd ,
		\end{equation*}
		and the partial derivative with respect to $\sigma$ is given by the
		same formula with $\phi$ in place of $\varphi$. By \Ftwo\ and Lemma
		\ref{lem:f}\ref{it:f-sign}, we have $wf'(w)+f(w)\ge2f(w)>0$ in $\RN$,
		so that the latter is negative. By the implicit function theorem,
		there is a function $\sigma(t)$ of class $\cC^1$ near $t=0$, with
		$\sigma(0)=0$, such that
		\begin{equation*}
			M_0\big(w+t\varphi+\sigma(t)\phi\big)=0,
			\qquad
			\sigma'(0)=-\frac{\int_{\RN}\big(wf'(w)+f(w)\big)\varphi\dd}
			{\int_{\RN}\big(wf'(w)+f(w)\big)\phi\dd}.
		\end{equation*}
		Since $w+t\varphi+\sigma(t)\phi\in\cM_0$, the function
		$t\mapsto\frac12\norm[2]{w+t\varphi+\sigma(t)\phi}^2$ has a minimum at
		$t=0$, that is,
		$\int_{\RN}w\varphi\dd+\sigma'(0)\int_{\RN}w\phi\dd=0$. Hence
		\begin{equation*}
			\int_{\RN}w\varphi\dd=\lambda\int_{\RN}\big(wf'(w)+f(w)\big)\varphi\dd
			\qquad\text{for every }\varphi\in\cC_c^\infty(\RN),
		\end{equation*}
		where $\lambda>0$ is the quotient of $\int_{\RN}w\phi\dd$ and
		$\int_{\RN}\big(wf'(w)+f(w)\big)\phi\dd$, which does not depend on
		$\varphi$. Since $w$ is continuous and positive, it follows that
		\begin{equation*}
			1=\lambda\left(f'(w(x))+\frac{f(w(x))}{w(x)}\right)
			\qquad\text{for every }x\in\RN.
		\end{equation*}
		By \Fone, the right-hand side tends to zero as $|x|\to\infty$, because
		$w(x)\to0$. This contradiction shows that $\cz<c_*$.
	\end{proof}
	
	\begin{proposition}\label{prop:divergence}
		Assume \Fone--\Ffour. Suppose that $N\in\{3,4\}$, or that
		\eqref{eq:critical-comparison} holds. Then
		$c_-(\mu)\to+\infty$ as $\mu\to0^+$.
	\end{proposition}
	\begin{proof}
		Suppose that the assertion fails. Since $c_-(\mu)$ is attained for
		every $\mu>0$, there are $\mu_n\to0^+$ and $u_n\in\cC(\mu_n)$ such that
		$(u_n)$ is bounded in $L^2(\RN)$.
		
		Let $N\in\{3,4\}$. As explained at the beginning of this section, we
		may assume \eqref{eq:positive-radial}. Let $w$ be given by Lemma
		\ref{lem:locallimit}. By Lemma
		\ref{lem:locallimit}\ref{it:ll-lower},
		\begin{equation*}
			\int_{\RN}w^2\dd\ge c\int_1^\infty r^{2(2-N)}r^{N-1}\,dr
			=c\int_1^\infty r^{3-N}\,dr=\infty ,
		\end{equation*}
		which contradicts Lemma \ref{lem:locallimit}\ref{it:ll-basic}.
		
		Let now \eqref{eq:critical-comparison} hold, and recall that $p=2^*$ in
		this case. Since $\norm[\infty]{\nabla u_n}<1$, we have
		$(1+t)u_n\in\D$ for $|t|$ small. We take $v=(1+t)u_n$ in
		\eqref{eq:critical}, divide by $t$ and let $t\to0^\pm$, which is
		justified by the dominated convergence theorem. This gives
		\begin{equation}\label{eq:endpoint-critical-nehari}
			\int_{\RN}\frac{|\nabla u_n|^2}{\sqrt{1-|\nabla u_n|^2}}\dd
			+\mu_n\norm[2]{u_n}^2=\int_{\RN}f(u_n)u_n\dd .
		\end{equation}
		On the other hand, $M_{\mu_n}(u_n)=0$ and
		\eqref{eq:critical-comparison} give
		\begin{equation*}
			2\Psi(u_n)+\frac{N+2}{N}\mu_n\norm[2]{u_n}^2
			=2\int_{\RN}F(u_n)\dd+\frac2N\int_{\RN}f(u_n)u_n\dd
			\ge\int_{\RN}f(u_n)u_n\dd .
		\end{equation*}
		Since
		\begin{equation*}
			\frac{t^2}{\sqrt{1-t^2}}-2\big(1-\sqrt{1-t^2}\big)
			=\frac{t^4}{\sqrt{1-t^2}\big(1+\sqrt{1-t^2}\big)^2}
			\ge\frac{t^4}{4}\qquad\text{for }0\le t<1,
		\end{equation*}
		the last inequality and \eqref{eq:endpoint-critical-nehari} imply
		\begin{equation*}
			\frac14\norm[4]{\nabla u_n}^4
			\le\int_{\RN}\frac{|\nabla u_n|^2}{\sqrt{1-|\nabla u_n|^2}}\dd
			-2\Psi(u_n)
			\le\frac{2}{N}\mu_n\norm[2]{u_n}^2\to0 .
		\end{equation*}
		By the Gagliardo--Nirenberg inequality, whose exponents are admissible
		since $\frac1{2^*}=\vartheta\big(\frac14-\frac1N\big)+\frac{1-\vartheta}{2}$
		for $\vartheta=\frac4{N+4}$,
		\begin{equation*}
			\norm[2^*]{u_n}\le C\norm[4]{\nabla u_n}^{\frac4{N+4}}
			\norm[2]{u_n}^{\frac N{N+4}}\to0 .
		\end{equation*}
		Since $(u_n)$ is bounded in $L^\infty(\RN)$ by Lemma
		\ref{lem:amplitude}, Lemma \ref{lem:f}\ref{it:f-local} with $q=2^*$,
		\eqref{eq:actionidentity} and Lemma \ref{prop:positive} give
		\begin{equation*}
			0<a_0\le a(\mu_n)\le\frac1N\int_{\RN}f(u_n)u_n\dd
			\le C\norm[2^*]{u_n}^{2^*}\to0,
		\end{equation*}
		a contradiction.
	\end{proof}
	
	\begin{lemma}\label{lem:decay}
		Let $A>0$ and $\alpha>0$ be such that $(p-2)\alpha>2$, where $p$ is the
		exponent in \Ffour. Then there is $C>0$, depending only on $N$, $f$,
		$A$ and $\alpha$, with the following property. If $\mu>0$ and
		$u(x)=W(|x|)\in\cC(\mu)$ is positive with
		\begin{equation*}
			W(r)\le A\min\{1,r^{-\alpha}\}\qquad\text{for }r>0,
		\end{equation*}
		then
		\begin{equation}\label{eq:decay}
			W(r)\le Cr^{2-N}\qquad\text{and}\qquad |W'(r)|\le Cr^{1-N}
			\qquad\text{for }r\ge1.
		\end{equation}
	\end{lemma}
	\begin{proof}
		In this proof the constants $C$ depend only on $N$, $f$, $A$, $\alpha$
		and on the exponent $\beta$ below. By Lemma
		\ref{lem:f}\ref{it:f-local} with $q=p$, $f(s)\le Cs^{p-1}$ for
		$0\le s\le A$. Let $\gamma=\frac12\big((p-2)\alpha-2\big)>0$, and
		suppose that
		\begin{equation}\label{eq:decay-step}
			W(r)\le Cr^{-\beta}\qquad\text{for }r\ge1
		\end{equation}
		with some $\beta\ge\alpha$. By \eqref{eq:fluxbound}, for $r\ge1$,
		\begin{equation*}
			-r^{N-1}W'(r)\le\int_0^rs^{N-1}f(W(s))\,ds
			\le C\left(1+\int_1^rs^{N-1-(p-1)\beta}\,ds\right).
		\end{equation*}
		If $(p-1)\beta>N$, then the right-hand side is bounded, so that
		$-W'(r)\le Cr^{1-N}$ for $r\ge1$, and integrating over $(r,\infty)$ we
		get \eqref{eq:decay}. If $(p-1)\beta\le N$, then
		\begin{equation*}
			\int_1^rs^{N-1-(p-1)\beta}\,ds
			\le\int_1^rs^{N-1-(p-1)\beta+\gamma}\,ds
			\le\frac1\gamma\,r^{N-(p-1)\beta+\gamma},
		\end{equation*}
		hence $-W'(r)\le Cr^{1-(p-1)\beta+\gamma}$ for $r\ge1$. Since
		\begin{equation*}
			(p-1)\beta-\gamma-2=\beta+(p-2)\beta-2-\gamma\ge\beta+\gamma,
		\end{equation*}
		an integration over $(r,\infty)$ yields \eqref{eq:decay-step} with
		$\beta+\gamma$ in place of $\beta$. We start from $\beta=\alpha$ and
		apply this argument with the exponents $\alpha+j\gamma$,
		$j=0,1,\dots$. After finitely many steps, whose number depends only on
		$N$, $p$ and $\alpha$, we arrive at the case $(p-1)\beta>N$, and
		\eqref{eq:decay} follows.
	\end{proof}
	
	\begin{proposition}\label{prop:bounded-mass-limit}
		Let $N\ge5$, $\mu_n\to0^+$ and $u_n\in\cC(\mu_n)$, and suppose that
		$(u_n)$ is bounded in $L^2(\RN)$. Then there are $y_n\in\RN$ and
		$w\in\cM_0$ such that, up to a subsequence, $u_n(\cdot+y_n)\to w$ in
		$H^1(\RN)$. Moreover, $w$ is a classical solution of
		\eqref{eq:zerofreq}, $\cE(w)=a_0$, $\norm[\infty]{\nabla w}<1$, and
		$w(x)=\pm W(|x|)$, where $W$ is positive and nonincreasing on
		$[0,\infty)$ and
		\begin{equation}\label{eq:w-decay}
			c\,r^{2-N}\le W(r)\le Cr^{2-N},\qquad|W'(r)|\le Cr^{1-N}
			\qquad\text{for }r\ge1
		\end{equation}
		with some constants $C,c>0$.
	\end{proposition}
	\begin{proof}
		As explained at the beginning of this section, it is
		sufficient to prove the assertion with $y_n=0$ under the assumption
		\eqref{eq:positive-radial}. Let $W$ and $w$ be given by Lemma
		\ref{lem:locallimit}. Then, up to a subsequence, $W_n\to W$ in
		$\cC^1([0,R])$ for every $R>0$, the function $w$ is a classical
		solution of \eqref{eq:zerofreq} with $\norm[\infty]{\nabla w}<1$, and
		$W$ is positive, nonincreasing and satisfies the lower bound in
		\eqref{eq:w-decay}. Since $W_n$ is nonincreasing and $(u_n)$ is bounded
		in $L^2(\RN)$ and in $L^\infty(\RN)$, we have
		$W_n(r)\le C\min\{1,r^{-N/2}\}$. Moreover,
		$(p-2)\frac N2\ge\frac{2N}{N-2}>2$ because $p\ge2^*$. Hence Lemma
		\ref{lem:decay} with $\alpha=N/2$ gives a constant $C$, independent of
		$n$, such that
		\begin{equation}\label{eq:uniformtails}
			W_n(r)\le Cr^{2-N},\qquad|W_n'(r)|\le Cr^{1-N}\qquad(r\ge1).
		\end{equation}
		Letting $n\to\infty$ we obtain the same bounds for $W$, which completes
		the proof of \eqref{eq:w-decay}. The functions $r^{2(2-N)}r^{N-1}$ and
		$r^{2(1-N)}r^{N-1}$ are integrable on $(1,\infty)$ for $N\ge5$. Hence
		\eqref{eq:uniformtails}, the convergence $W_n\to W$ in $\cC^1([0,R])$
		for every $R$ and the dominated convergence theorem give $u_n\to w$
		and $\nabla u_n\to\nabla w$ in $L^2(\RN)$. In particular,
		$w\in\D\setminus\{0\}$. Since $0\le1-\sqrt{1-|\xi|^2}\le|\xi|^2$ for
		$|\xi|\le1$, the dominated convergence theorem also gives
		$\Psi(u_n)\to\Psi(w)$. Furthermore, by Lemma
		\ref{lem:f}\ref{it:f-local} with $q=2^*$ and \eqref{eq:uniformtails},
		$F(u_n)+f(u_n)u_n\le C|x|^{-2N}$ for $|x|\ge1$, so that
		\begin{equation*}
			\int_{\RN}F(u_n)\dd\to\int_{\RN}F(w)\dd,\qquad
			\int_{\RN}G(u_n)\dd\to\int_{\RN}G(w)\dd .
		\end{equation*}
		Since $\mu_n\norm[2]{u_n}^2\to0$, passing to the limit in
		$M_{\mu_n}(u_n)=0$ and in $I_{\mu_n}(u_n)=a(\mu_n)$ we obtain
		$M_0(w)=0$ and $\cE(w)=\lim_{n\to\infty}a(\mu_n)=a_0$. Thus
		$w\in\cM_0$.
	\end{proof}

	In the rest of this section we write
	\begin{equation*}
		\ell_0:=\liminf_{\mu\to0^+}c_-(\mu)\in[c_*,+\infty].
	\end{equation*}
	
	\begin{lemma}\label{lem:frequency-gap}
		Let $N\in\{5,6\}$ and suppose that $\ell_0<\infty$. Then there is
		$\mu>0$ such that $a(\mu)<a_0+\ell_0\mu$.
	\end{lemma}
	\begin{proof}
		Choose $\mu_n\to0^+$ and $u_n\in\cC(\mu_n)$ with
		$\frac12\norm[2]{u_n}^2=c_-(\mu_n)\to\ell_0$, and let $w$ be given by
		Proposition \ref{prop:bounded-mass-limit}. Then
		$\frac12\norm[2]{w}^2=\ell_0$. Replacing $f$ by $s\mapsto-f(-s)$ and
		$w$ by $-w$ if necessary, we may assume that $w(x)=W(|x|)$, see the
		beginning of this section. Since $M_0(w)=0$, we have
		$\int_{\RN}f(w)w\dd=N\big(\Psi(w)-\int_{\RN}F(w)\dd\big)=Na_0$. In this
		proof the constants $C$ and $c$ do not depend on $\eps$.
		
		\emph{Step 1.} For $\eps\in(0,\frac1{16}]$ let
		$r_\eps=(\eps^{-2}+\eps^{-1})^{1/2}$ and
		\begin{equation*}
			\eta(x)=\min\big\{\eps^2\max\{|x|^2-\eps^{-1},0\},1\big\},
			\qquad h=\eta w .
		\end{equation*}
		Thus $\eta=0$ for $|x|\le\eps^{-1/2}$, $\eta=1$ for $|x|\ge r_\eps$,
		$|\nabla\eta|\le2\eps^2|x|$, and $\eps^{-1}\le r_\eps\le2\eps^{-1}$.
		Consider
		\begin{equation*}
			\delta_\eps=\int_{\RN}\eta w^2\dd,\qquad
			d_\eps=\int_{\RN}|\nabla h|^2\dd,\qquad
			t_\eps=\int_{|x|\ge\eps^{-1/2}}\big(F(w)+f(w)w\big)\dd .
		\end{equation*}
		We claim that, as $\eps\to0^+$,
		\begin{equation}\label{eq:smallness}
			\delta_\eps\to0,\qquad\frac{d_\eps}{\delta_\eps^2}\to0,\qquad
			\frac{t_\eps}{\delta_\eps}\to0 .
		\end{equation}
		Indeed, $\delta_\eps\le\int_{|x|\ge\eps^{-1/2}}w^2\dd\to0$, and Lemma
		\ref{lem:f}\ref{it:f-local} and \eqref{eq:w-decay} give
		$F(w)+f(w)w\le Cw^{2^*}\le C|x|^{-2N}$ for $|x|\ge1$, hence
		$t_\eps\le C\eps^{N/2}$. Since
		$|\nabla h|^2\le2\eta^2|\nabla w|^2+2|\nabla\eta|^2w^2$ and
		$\eta\le\eps^2|x|^2$, we obtain from \eqref{eq:w-decay}
		\begin{equation*}
			d_\eps\le C\eps^4\int_{\eps^{-1/2}}^{r_\eps}r^{6-2N}r^{N-1}\,dr
			+C\int_{r_\eps}^\infty r^{2-2N}r^{N-1}\,dr
			\le\begin{cases}C\eps^3,&N=5,\\ C\eps^4|\log\eps|,&N=6.\end{cases}
		\end{equation*}
		On the other hand, if $N=5$, then
		\begin{equation*}
			\delta_\eps\ge\int_{|x|\ge r_\eps}w^2\dd
			\ge c\int_{r_\eps}^\infty r^{-2}\,dr
			=\frac{c}{r_\eps}\ge c\,\eps .
		\end{equation*}
		If $N=6$, then $\eta\ge\frac34\eps^2|x|^2$ for
		$2\eps^{-1/2}\le|x|\le\eps^{-1}$, and
		\begin{equation*}
			\delta_\eps\ge c\,\eps^2\int_{2\eps^{-1/2}}^{\eps^{-1}}\frac{dr}{r}
			\ge c\,\eps^2|\log\eps| .
		\end{equation*}
		These estimates give \eqref{eq:smallness}.
		
		\emph{Step 2.} Let $z=w-h=(1-\eta)w$. Then $0\le z\le w$, $z=w$ for
		$|x|\le\eps^{-1/2}$ and $z=0$ for $|x|\ge r_\eps$. Moreover,
		$|\nabla z|=|\nabla w|$ for $|x|<\eps^{-1/2}$, and by
		\eqref{eq:w-decay}, for $|x|\ge\eps^{-1/2}$,
		\begin{equation*}
			|\nabla z|\le|\nabla w|+|\nabla\eta|w
			\le C\eps^{(N-1)/2}+C\eps^{(N+1)/2} .
		\end{equation*}
		Hence, for small $\eps$, $z$ is a Lipschitz function with compact
		support and $\norm[\infty]{\nabla z}\le\norm[\infty]{\nabla w}<1$. In
		particular, $z\in\D\setminus\{0\}$. The Hessian of
		$\xi\mapsto1-\sqrt{1-|\xi|^2}$ is bounded on the ball
		$\{|\xi|\le\norm[\infty]{\nabla w}\}$, so that Taylor's formula gives
		\begin{equation*}
			1-\sqrt{1-|\zeta|^2}\le1-\sqrt{1-|\xi|^2}
			+\frac{\xi\cdot(\zeta-\xi)}{\sqrt{1-|\xi|^2}}
			+C|\zeta-\xi|^2\qquad\text{for }|\xi|,|\zeta|\le\norm[\infty]{\nabla w} .
		\end{equation*}
		We apply this inequality with $\xi=\nabla w$ and
		$\zeta=\nabla z=\nabla w-\nabla h$ and integrate. Multiplying
		\eqref{eq:zerofreq} by the bounded Lipschitz function $h$ and
		integrating by parts over the ball of radius $R$, where the boundary
		term is $O(R^{2-N})$ by \eqref{eq:w-decay}, we get
		$\int_{\RN}\nabla w\cdot\nabla h\,(1-|\nabla w|^2)^{-1/2}\dd
		=\int_{\RN}f(w)h\dd$. Therefore
		\begin{equation}\label{eq:Psi-v1}
			\Psi(z)\le\Psi(w)-\int_{\RN}f(w)h\dd+Cd_\eps
			\le\Psi(w)+Cd_\eps .
		\end{equation}
		Since $F$ is convex, we have $F(z)\ge F(w)-f(w)h$, and
		\eqref{eq:Psi-v1} gives
		\begin{equation}\label{eq:E-v1}
			\cE(z)\le a_0+Cd_\eps .
		\end{equation}
		Moreover, since $z=w$ for $|x|\le\eps^{-1/2}$,
		$\Psi(w)=\int_{\RN}G(w)\dd$ and $(1-\eta)^2\le1-\eta$,
		\begin{equation}\label{eq:v1-integrals}
			\int_{\RN}f(z)z\dd\ge Na_0-t_\eps,\qquad
			\int_{\RN}G(z)\dd\ge\Psi(w)-t_\eps,\qquad
			\norm[2]{z}^2\le\norm[2]{w}^2-\delta_\eps .
		\end{equation}
		
		\emph{Step 3.} Let
		\begin{equation*}
			\tau=\frac{\delta_\eps}{(N+2)2^{N+2}\norm[2]{w}^2},\qquad
			\theta=1+\tau,\qquad v=z_\theta ,
		\end{equation*}
		and let $\eps$ be so small that $\tau\le1$. By the convexity of $F$,
		$F(\theta s)\ge F(s)+\tau f(s)s$ for $s\ge0$. Hence, by
		\eqref{eq:scaling}, \eqref{eq:E-v1}, \eqref{eq:v1-integrals} and the
		elementary inequality $(1+\tau)^N(1-N\tau)\le1$,
		\begin{equation}\label{eq:E-v}
			\cE(v)=\theta^N\Big(\Psi(z)-\int_{\RN}F(\theta z)\dd\Big)
			\le\theta^N\big(a_0(1-N\tau)+Cd_\eps+\tau t_\eps\big)
			\le a_0+C\big(d_\eps+\delta_\eps t_\eps\big).
		\end{equation}
		By Lemma \ref{lem:f}\ref{it:f-mono},
		$G(\theta s)\ge\theta^2G(s)\ge(1+2\tau)G(s)$, and by
		\eqref{eq:Psi-v1}, \eqref{eq:v1-integrals} and \eqref{eq:smallness},
		for small $\eps$,
		\begin{equation*}
			\int_{\RN}G(v)\dd-\Psi(v)
			=\theta^N\Big(\int_{\RN}G(\theta z)\dd-\Psi(z)\Big)
			\ge2\tau\Psi(w)-3t_\eps-Cd_\eps\ge\tau\Psi(w).
		\end{equation*}
		Since
		$\norm[2]{v}^2=\theta^{N+2}\norm[2]{z}^2\le2^{N+2}\norm[2]{w}^2$, the
		number
		\begin{equation*}
			\mu_v:=\frac{2N}{N+2}\cdot
			\frac{\int_{\RN}G(v)\dd-\Psi(v)}{\norm[2]{v}^2}
		\end{equation*}
		satisfies $\mu_v\ge c\,\delta_\eps$, and $M_{\mu_v}(v)=0$ by
		\eqref{eq:M}, that is, $v\in\cM_{\mu_v}$. Furthermore, since
		$\theta^{N+2}-1\le(N+2)2^{N+1}\tau$, by \eqref{eq:v1-integrals} and the
		choice of $\tau$,
		\begin{equation*}
			\norm[2]{v}^2=\theta^{N+2}\norm[2]{z}^2
			\le\norm[2]{w}^2-\delta_\eps+(\theta^{N+2}-1)\norm[2]{w}^2
			\le\norm[2]{w}^2-\frac12\delta_\eps .
		\end{equation*}
		Combining this inequality with \eqref{eq:E-v} we conclude that
		\begin{equation*}
			a(\mu_v)\le I_{\mu_v}(v)
			=\cE(v)+\frac{\mu_v}{2}\norm[2]{v}^2
			\le a_0+\frac{\mu_v}{2}\norm[2]{w}^2
			+C\big(d_\eps+\delta_\eps t_\eps\big)-c\,\delta_\eps^2 .
		\end{equation*}
		By \eqref{eq:smallness}, the sum of the last two terms is negative for
		small $\eps$. Since $\frac12\norm[2]{w}^2=\ell_0$, the lemma follows
		with $\mu=\mu_v$.
	\end{proof}
	
	\begin{lemma}\label{prop:cstar-rho}
		Assume \Fone--\Ffour.
		\begin{enumerate}[label=\textnormal{(\alph*)}]
			\item\label{it:cr-strict}
			If $c_*<\ell_0$, then the infimum defining $c_*$ is attained and
			$c_*<\rho_0$.
			\item\label{it:cr-rho}
			If $\rho_0<\ell_0$, then the infimum defining $\rho_0$ is attained.
		\end{enumerate}
	\end{lemma}
	\begin{proof}
		By \eqref{eq:integral0}, the number $\rho(\mu)$ in \eqref{eq:q} is the
		mean value of $c_-$ over $(0,\mu)$. Hence $\rho$ is continuous,
		$\liminf_{\mu\to0^+}\rho(\mu)\ge\ell_0$, and
		$\rho(\mu)\to+\infty$ as $\mu\to\infty$ by Proposition
		\ref{prop:fixed}\ref{it:large}.
		
		\ref{it:cr-strict} Let $c_*<\ell<\ell_0$. By the definition of $\ell_0$
		and Proposition \ref{prop:fixed}\ref{it:large}, there are
		$0<\mu_1<\mu_2$ such that $c_-\ge\ell$ on
		$(0,\mu_1)\cup(\mu_2,\infty)$. Hence $c_*$ is the infimum of $c_-$
		on $[\mu_1,\mu_2]$, and it is attained since $c_-$ is lower
		semicontinuous by Proposition \ref{prop:fixed}\ref{it:compact}.
		Moreover, $\rho(\mu)\ge\ell$ for $\mu\le\mu_1$, and for $\mu>\mu_1$
		\begin{equation*}
			\rho(\mu)=\frac1\mu\int_0^\mu c_-(t)\,dt
			\ge\frac{\ell\mu_1+c_*(\mu-\mu_1)}{\mu}
			=c_*+\frac{(\ell-c_*)\mu_1}{\mu} .
		\end{equation*}
		Since $\rho(\mu)\to+\infty$ as $\mu\to\infty$, these bounds give
		$\rho_0>c_*$.
		
		\ref{it:cr-rho} Since $\rho$ is continuous,
		$\liminf_{\mu\to0^+}\rho(\mu)\ge\ell_0>\rho_0$ and
		$\rho(\mu)\to+\infty$ as $\mu\to\infty$, the infimum of $\rho$ is
		attained.
	\end{proof}
	
	\begin{proof}[Proof of Theorem \ref{thm:attainment}]
		\ref{it:attain} Let $3\le N\le6$, or suppose that
		\eqref{eq:critical-comparison} holds. We claim that $\rho_0<\ell_0$.
		This is clear if $\ell_0=+\infty$. If $\ell_0<\infty$, then
		$N\in\{5,6\}$ by Proposition \ref{prop:divergence}, and Lemma
		\ref{lem:frequency-gap} gives $\rho(\mu)<\ell_0$ for some $\mu>0$. Since
		$c_*\le\rho_0$ by Lemma \ref{prop:thresholds}, we also have
		$c_*<\ell_0$. Hence, by Lemma
		\ref{prop:cstar-rho}\ref{it:cr-strict} and \ref{it:cr-rho}, the infima
		defining $c_*$ and $\rho_0$ are attained, and $c_*<\rho_0$.
		
		\ref{it:notattain} Suppose that the infimum defining $c_*$ is not
		attained. Then $N\ge7$ by \ref{it:attain}, and $\ell_0=c_*$, because
		$\ell_0\ge c_*$ and the strict inequality would give the attainment by
		Lemma \ref{prop:cstar-rho}\ref{it:cr-strict}. Choose
		$\mu_n\to0^+$ and $u_n\in\cC(\mu_n)$ with
		$\frac12\norm[2]{u_n}^2=c_-(\mu_n)\to c_*$. By Proposition
		\ref{prop:bounded-mass-limit}, there are $y_n\in\RN$
		and $w\in\cM_0$ such that, up to a subsequence, $u_n(\cdot+y_n)\to w$
		in $H^1(\RN)$, where $w$ is a radial classical solution of
		\eqref{eq:zerofreq} of constant sign with $\cE(w)=a_0$. In particular,
		$\frac12\norm[2]{w}^2=c_*$.
	\end{proof}
	
	\appendix
	\section{A multiplier rule for constrained local minima}\label{app:multiplier}
	
	In this appendix $X$ denotes an arbitrary real Banach space and $\Psi$
	an arbitrary convex functional on $X$. The space $X$ and the functional
	$\Psi$ of Section \ref{sec:intro} are a particular case. We consider
	local minimum points of $\Psi-\Phi$ under a constraint $Q=q$, where
	$\Phi$ and $Q$ are functionals of class $\cC^1$. Since $\Psi$ need not
	be differentiable, the Lagrange multiplier rule is not available in
	general. The following theorem gives a multiplier, together with its
	value, if the problem admits a scaling which acts on $\Psi$ by a
	multiplicative factor and changes the constraint at first order.
	
	\begin{theorem}\label{thm:multiplier}
		Let $X$ be a real Banach space, let $\Psi:X\to(-\infty,+\infty]$ be a
		convex functional with the effective domain
		$D(\Psi)=\{w\in X:\Psi(w)<+\infty\}$, and let $\Phi,Q\in\cC^1(X,\R)$.
		Let $J\subset(0,\infty)$ be an open interval containing $1$, and let
		$(T_\theta)_{\theta\in J}$ be bounded linear operators on $X$ such that
		$T_1$ is the identity and the following conditions hold.
		\begin{enumerate}[label=\textnormal{(\alph*)}]
			\item\label{it:T-cont}
			$\sup_{\theta\in J}\norm[\mathcal L(X)]{T_\theta}<\infty$, and
			$T_\theta w\to w$ in $X$ as $\theta\to1$ for every $w\in X$.
			\item\label{it:T-Psi}
			$T_\theta w\in D(\Psi)$ and $\Psi(T_\theta w)=\gamma(\theta)\Psi(w)$
			for every $w\in D(\Psi)$ and $\theta\in J$, where
			$\gamma:J\to(0,\infty)$ is differentiable at $1$ and $\gamma(1)=1$.
			\item\label{it:T-C1}
			The functions $(\theta,w)\mapsto\Phi(T_\theta w)$ and
			$(\theta,w)\mapsto Q(T_\theta w)$ are of class $\cC^1$ on $J\times X$.
		\end{enumerate}
		Let $q\in\R$, and let $u\in D(\Psi)$ with $Q(u)=q$ be a local minimum
		point of $\Psi-\Phi$ on the set $\{w\in D(\Psi):Q(w)=q\}$ in the
		topology of $X$. If
		$\frac{d}{d\theta}Q(T_\theta u)\big|_{\theta=1}\ne0$, then
		\begin{equation}\label{eq:multiplier-rule}
			\Psi(v)-\Psi(u)-\Phi'(u)[v-u]+\lambda\,Q'(u)[v-u]\ge0
			\qquad\text{for every }v\in X,
		\end{equation}
		where
		\begin{equation}\label{eq:multiplier}
			\lambda=\frac{\frac{d}{d\theta}\Phi(T_\theta u)\big|_{\theta=1}
				-\gamma'(1)\Psi(u)}
			{\frac{d}{d\theta}Q(T_\theta u)\big|_{\theta=1}}.
		\end{equation}
	\end{theorem}
	\begin{proof}
		If $v\notin D(\Psi)$, then \eqref{eq:multiplier-rule} is trivial. Let
		$v\in D(\Psi)$, and for $t\in\R$ put $z_t=u+t(v-u)$. Since $D(\Psi)$ is
		convex, $z_t\in D(\Psi)$ for $0\le t\le1$. By \ref{it:T-C1}, the
		function $(\theta,t)\mapsto Q(T_\theta z_t)$ is of class $\cC^1$ on
		$J\times\R$. It equals $q$ at $(1,0)$, and its partial derivatives at
		this point are $\frac{d}{d\theta}Q(T_\theta u)\big|_{\theta=1}\ne0$ and
		$Q'(u)[v-u]$. By the implicit function theorem, there are
		$t_0\in(0,1]$ and a function $\theta:(-t_0,t_0)\to J$ of class $\cC^1$
		such that $\theta(0)=1$,
		\begin{equation}\label{eq:multiplier-theta}
			Q\big(T_{\theta(t)}z_t\big)=q\quad\text{for }|t|<t_0,\qquad
			\theta'(0)=-\frac{Q'(u)[v-u]}
			{\frac{d}{d\theta}Q(T_\theta u)\big|_{\theta=1}} .
		\end{equation}
		For $0\le t<t_0$ let $w_t=T_{\theta(t)}z_t$. Then $w_t\in D(\Psi)$ by
		\ref{it:T-Psi}, $Q(w_t)=q$, and by \ref{it:T-cont}
		\begin{equation*}
			\norm[X]{w_t-u}\le\norm[\mathcal L(X)]{T_{\theta(t)}}\norm[X]{z_t-u}
			+\norm[X]{T_{\theta(t)}u-u}\to0\qquad\text{as }t\to0^+ .
		\end{equation*}
		Hence, for small $t\ge0$, the local minimality of $u$, \ref{it:T-Psi}
		and the convexity of $\Psi$ give
		\begin{equation*}
			\Psi(u)-\Phi(u)\le\Psi(w_t)-\Phi(w_t)
			\le\gamma(\theta(t))\big((1-t)\Psi(u)+t\Psi(v)\big)
			-\Phi\big(T_{\theta(t)}z_t\big).
		\end{equation*}
		By \ref{it:T-Psi} and \ref{it:T-C1}, the right-hand side is
		differentiable at $t=0$, where it equals $\Psi(u)-\Phi(u)$. Therefore
		its derivative at $t=0$ is nonnegative, that is,
		\begin{equation*}
			\Psi(v)-\Psi(u)-\Phi'(u)[v-u]
			+\theta'(0)\Big(\gamma'(1)\Psi(u)
			-\frac{d}{d\theta}\Phi(T_\theta u)\Big|_{\theta=1}\Big)\ge0 .
		\end{equation*}
		Inserting the value of $\theta'(0)$ from \eqref{eq:multiplier-theta}
		we obtain \eqref{eq:multiplier-rule} with $\lambda$ given by
		\eqref{eq:multiplier}.
	\end{proof}
	
	Inequality \eqref{eq:multiplier-rule} means that $u$ is a critical
	point, in the sense of \cite{S}, of the functional
	$\Psi-\Phi+\lambda Q$.

\medskip
\noindent\textbf{Acknowledgements.} J. M. was partly supported by the National Science Centre, Poland, under the research grant 2025/57/B/ST1/04804.

\smallskip


	\smallskip
	\noindent\textbf{Data availability.} The manuscript has no associated data.
	
	\smallskip
	\noindent\textbf{Conflict of interest.} The authors declare no conflict of interest.

\end{document}